\documentclass[reqno,final]{amsart}
\usepackage{natbib}  
\usepackage{fancyhdr} 
\usepackage{color} 
\usepackage{hyperref} 
\usepackage{graphicx} 

\usepackage{pstricks}
\usepackage{amssymb}

\definecolor{aleacolor}{rgb}{0.16,0.59,0.78}

\hypersetup{
breaklinks,
colorlinks=true,
linkcolor=aleacolor,
urlcolor=aleacolor,
citecolor=aleacolor}

\renewcommand{\cite}{\citet}
\renewcommand{\cite}{\citep}

\theoremstyle{plain}
\newtheorem{theorem}{Theorem}[section]                                          
\newtheorem{proposition}[theorem]{Proposition}                          

\newtheorem{corollary}[theorem]{Corollary}

\theoremstyle{definition}
\newtheorem{definition}[theorem]{Definition}
\theoremstyle{remark}
\newtheorem{remark}[theorem]{Remark}

\newtheorem{cor}{Corollary}
\newtheorem{theo}{Theorem}

\newtheorem{prop}{Proposition}

\newtheorem{rem}{Remark}
\newtheorem{pty}{Property}
\newtheorem{def-prop}{Definition-Property}

\makeatletter \@addtoreset{equation}{section} \makeatother

\begin{document}

\title[A third-moment theorem]{A third-moment theorem and precise asymptotics for variations of stationary Gaussian sequences}

\author{L\'{e}o Neufcourt}
\author{Frederi G. Viens}

\address{Department of Statistics, Columbia University\newline
1255 Amsterdam Avenue,\newline
Room 1005 SSW, MC 4690\newline
New York, NY, USA.}

\address{Ex-affiliate of CMAP, Ecole Polytechnique\newline
Route de Saclay,\newline
UMR CNRS 7641\newline
Palaiseau, France.}

\address{Department of Statistics, Purdue University\newline
150 N. University St.\newline
Office MATH 200\newline
West Lafayette, IN, USA.}

\email{leo.neufcourt@columbia.edu, viens@stat.purdue.edu}
\urladdr{\url{http://www.stat.purdue.edu/~viens/}}

\thanks{Research supported by CIMFAV at Universidad de Valparaiso, Ecole Polytechnique}

\subjclass[2010 MSC]{60G15, 60F05, 60H07, 60G22.} 
\keywords{Stationary Gaussian Process, Wiener Space, Central Limit Theorem, Berry-Ess\'{e}en, Breuer-Major, Second Chaos.}

\begin{abstract}
In two new papers \cite{alea} and \cite{4M}, sharp general quantitative
bounds \ are given to complement the well-known fourth moment theorem of
Nualart and Peccati, by which a sequence in a fixed Wiener chaos converges
to a normal law if and only if its fourth cumulant converges to $0$. The
bounds show that the speed of convergence is precisely of order the maximum
of the fourth cumulant and the absolute value of the third moment
(cumulant). Specializing to the case of normalized centered quadratic
variations for stationary Gaussian sequences, we show that a third moment
theorem holds: convergence occurs if and only if the sequence's third
moments tend to $0$. This is proved for sequences with general decreasing
covariance, by using the result of \cite{4M}, and finding the exact speed of
convergence to $0$ of the quadratic variation's third and fourth cumulants. 
\cite{4M} also allows us to derive quantitative estimates for the speeds of
convergence in a class of log-modulated covariance structures, which puts in
perspective the notion of critical Hurst parameter when studying the
convergence of fractional Brownian motion's quadratic variation. We also
study the speed of convergence when the limit is not Gaussian but rather a
second-Wiener-chaos law. Using a log-modulated class of spectral densities,
we recover a classical result of Dobrushin-Major/Taqqu whereby the limit is
a Rosenblatt law, and we provide new convergence speeds. The conclusion in
this case is that the price to pay to obtain a Rosenblatt limit despite a
slowly varying modulation is a very slow convergence speed, roughly of the
same order as the modulation.
\end{abstract}

\maketitle

\newpage

\section{Introduction}

We are inspired by the following reformulation of Theorem 1.2 in \cite{4M},
which is itself based on ideas contained in \cite{alea}.

\begin{theorem}[$4^{th}$ moment theorem in total variation and convergence rates%
]
\label{T4M}

If $(F_{n})_{n\geq 0}$ is a sequence in a fixed Wiener chaos (e.g. in the
second chaos), and $Var\left[ F_{n}\right] =1$, then $(F_{n})_{n\geq 0}$
converges in law towards $\mathcal{N}(0,1)$ if and only if $\mathbf{E}%
[F_{n}^{4}]\rightarrow 3=\mathbf{E}[N^4]$, where $N\sim \mathcal{N}(0,1)$.
Moreover the convergence rate in this case is $M_{n}:=\max (\mathbf{E}%
[F_{n}^{4}]-3,\left\vert \mathbf{E}[F_{n}^{3}]\right\vert )$, in the sense
of commensurability for the total variation metric $d_{TV}(F_{n},N)\asymp
M_{n}$, i.e. 
\begin{equation}
\exists c,C>0:cM_{n}\leq d_{TV}(F_{n},N)\leq CM_{n}.
\label{convergence_exacte}
\end{equation}
\end{theorem}

The first part of this theorem is known as the 4th moment theorem, proved
originally by Nualart and Peccati in \cite{NuP}. The second part, i.e.
relation (\ref{convergence_exacte}), suggests that the third moment is just
as important as the 4th moment when investigating the normal convergence of
sequences in a fixed Wiener chaos. Theorem \ref{T4M} also provides new
information about a rather successful estimate for evaluating normal
convergence speeds, which was established thanks to a research program
started in 2008 by Nourdin and Peccati in \cite{NP}. Specifically, in \cite%
{NPRAoP} an upper bound on $d_{TV}(F_{n},N)$ of the form $C\sqrt{\mathbf{E}%
[F_{n}^{4}]-3}$ is established (see Theorem 5.2.6 in \cite{npBook}); Theorem %
\ref{T4M} above shows that this estimate is not sharp in cases where the
third moment is dominated by the 4th moment minus 3 (a.k.a. the 4th cumulant 
$\kappa _{4}\left( F_{n}\right) $), and in the other cases, leaves the
question of sharpness of past results to a comparison with the third moment.
Such a discovery begs the question of how much one might improve certain
convergence results, e.g. of Berry-Ess\'{e}en type, by using Theorem \ref%
{T4M} instead of Theorem 5.2.6 in \cite{npBook}. The authors of \cite{npBook}
joined forces with Bierm\'{e} and Bonami to produce the first positive
result in this direction, in \cite{alea}: they worked with a weaker notion
of convergence than total-variation convergence, but were able to show, in
the case of the power variations of discrete fractional Brownian motion,
that the third moment seems to dominate the 4th cumulant in many cases, and
therefore determines normal convergence in those cases, yielding much
improved speeds as a consequence. In \cite{4M}, Nourdin and Peccati improve
this result by showing that it holds for the total variation distance,
thanks to the general Theorem \ref{T4M}.

Herein, we too base our analysis on the power afforded by Theorem \ref{T4M},
and focus our attention on broadening the study of quadratic variations from 
\cite{alea}, to include general stationary Gaussian sequences with no
reference to H\"{o}lder-continuity or self-similarity.

In addition to what is described above, our motivation is to show that, in
the case of quadratic variations (which live in the 2nd chaos) the fourth
moment theorem can be replaced by a third moment theorem, with corresponding
quantitative estimates of the total variation distance and of the relation
between the 3rd and 4th moments. Another motivation is to keep our analysis
as general as possible within the confines of variations of stationary
Gaussian sequences, with as little assumptions about their covariance
structure as we can.

In other words, we consider a sequence of centered identically distributed
Gaussian random variables $\left( X_{n}\right) _{n\in \mathbf{Z}}$ for which
there exists a function $\rho $ on $\mathbf{Z}$ such that for all integers $%
k,n$, $\mathbf{E}\left[ X_{n}X_{n+k}\right] =\rho \left( k\right) $. We only
assume that $\rho $ is of constant sign, and $\left\vert \rho \right\vert $
decreases near $+\infty $. Note that $\rho $ is necessarily symmetric (even)
and of positive type (meaning that $\left( \rho \left( k-\ell \right)
\right) _{k,\ell }$ is a non-negative definite matrix). Without loss of
generality, we assume that $Var\left[ X_{n}\right] =\rho \left( 0\right) =1$
throughout. We define the normalized centered \emph{quadratic variation} 
\begin{equation*}
F_{n}:=\frac{V_{n}}{\sqrt{v_{n}}},
\end{equation*}%
where $V_{n}:=\frac{1}{\sqrt{n}}\sum_{k=0}^{n-1}(X_{k}^{2}-1)$, and $%
v_{n}:=E[V_{n}^{2}]$. We prove the equivalence of third and fourth moment
theorems for the limit of $F_{n}$ under this general framework (Theorem \ref%
{T3M} on page \pageref{T3M}).

The presence of the \textquotedblleft normalizing\textquotedblright\ term $1/%
\sqrt{n}$ in the definition of $V_{n}$ is not needed (Theorem \ref{T3M}
remains valid without it); it is a convention that comes from working with
sequences where $v_{n}$ is bounded, so that the normalization of the
discrete centered quadratic variation $\sum_{k=0}^{n-1}(X_{k}^{2}-1)$ needed
to obtain a normal limit is $1/\sqrt{n}$, as one would hope for in a
straightforward generalization of the central limit theorem. This situation
is the well-known framework for the classical Breuer-Major central limit
theorem (see \cite{BM} or Theorem 7.2.4 in \cite{npBook}). Some might argue
that the interesting cases are those for which $v_{n}$ tends to infinity.
See step 2 in the proof of Proposition \ref{logprop} in Section \ref{LOG}
for a class of log-modulated models where the transition from bounded to
unbounded $v_{n}$ occurs, i.e. where the Breuer-Major statement becomes
invalid. Generally speaking, we argue that this type of \textquotedblleft
phase transition\textquotedblright\ is an artefact of the model class one
uses.

We investigate this question of phase transitions\ in the Berry-Ess\'{e}en
rates for the Breuer-Major central limit theorem, i.e. the rate at which $%
d_{TV}(F_{n},N)$ converges to $0$. In the case of fractional Gaussian noise
(fGn, for which $\rho \left( k\right) $ is equivalent to $H\left(
2H-1\right) \left\vert k\right\vert ^{2H-2}$ for large $\left\vert
k\right\vert $), up until very recently, the rate of convergence of $%
d_{TV}(F_{n},N)$ to $0$ was known to be the classical Berry-Ess\'{e}en rate $%
1/\sqrt{n}$ only for $H<5/8$. Thanks to the results in \cite{4M} and \cite%
{alea}, it is now known that this rate hold up to $H<2/3$. The optimal rates
for fBm are even known, as can be seen in \cite{alea} (also see Proposition
4.3 in \cite{4M}): the rate $1/\sqrt{n}$ is optimal for $H<2/3$; then for $%
H=2/3$, the optimal rate is $n^{-1/2}\log ^{2}n$, whereas for $H\in (2/3,3/4)
$, the optimal rate is $n^{6H-9/2}$. Therefore, speaking strictly about the
fGn scale, this can be interpreted as a \textquotedblleft phase
transition\textquotedblright , or a critical threshold $H=2/3$. We argue
that for general sequences, insofar as optimal rates are given by (\ref%
{convergence_exacte}) in Theorem \ref{T4M} if they can be computed, the
notion of critical threshold is model-class-dependent.

Herein, we show that the so-called critical threshold above can be
investigated in more detail, to reveal a range of possibilities for the
convergence rate at the Berry-Ess\'{e}en critical threshold\ $H=2/3$, for a
class of covariance functions with log-modulation. Our tools also enable us
to compute rates for situations which were left out in \cite{4M} and \cite%
{alea}: a second critical threshold\ at $H=3/4$ reveals another range of
possibilities under our log-modulated class where normal convergence holds.
This subclass of models with \textquotedblleft $H=3/4$\textquotedblright\
contains models for which the Breuer-Major theorem holds, and others for
which it does not. When normal convergence holds but the classical
Breuer-Major normalization fails, we find very slow rates of convergence:
rate $\log ^{-\frac{3}{2}}(n)$ in all cases except one exceptional model
which benefits from a $\log (\log (n))^{-\frac{3}{2}}$ correction. For fGn,
which is a special case of our results, we provide a slight improvement on
the only result we are aware of (Corollary 7.4.3 in \cite{npBook}) which
deals with the case \textquotedblleft $H=3/4$\textquotedblright : scaled
quadratic variations for fGn with $H=3/4$ have a speed of normal convergence
of order $\log ^{-1}\left( n\right) $. All our results mentioned above are
sharp, and are given in Proposition \ref{logprop} on page \pageref{logprop}.

Both this class of examples, given in Section \ref{LOG}, and our generic
third moment theorem given in Section \ref{THM}, are based on a precise
asymptotic expression for the third moment of $F_{n}$, and a precise upper
bound for the fourth moment of $F_{n}$, both given for general $\rho ,$ in
Theorem \ref{theo-est} in Section \ref{CUM}. This intrinsic study shows that
there is nothing special about the power rate of correlation decay when
investigating central limit theorems for quadratic variations. We also note
that we are able to develop a general third-moment theorem without needing a
direct estimation of the fourth cumulant $\kappa _{4}\left( F_{n}\right) $,
giving more credence to the claim that the third moment theorem is the right
tool for studying centered quadratic variations.

The key to this surprising shortcut is Proposition \ref{PropKey}, which
shows that one can bound the third moment $\kappa _{3}\left( F_{n}\right) $
below by $n^{1/4}\kappa _{4}\left( F_{n}\right) ^{3/4}$. This implies that $%
\kappa _{3}\left( F_{n}\right) $ is always the dominant term in the exact
rate of convergence of $d_{TV}(F_{n},N)$ to zero, showing in particular that
if $\kappa _{3}\left( F_{n}\right) \rightarrow 0$ then $\kappa _{4}\left(
F_{n}\right) \rightarrow 0$, which is sufficient to conclude that $%
F_{n}\rightarrow N$. This in turn implies the exact total-variation
convergence rate (\ref{convergence_exacte}) from Theorem \ref{T4M}. The
converse statement, that $\kappa _{4}\left( F_{n}\right) \rightarrow 0$
implies $\kappa _{3}\left( F_{n}\right) \rightarrow 0$, follows immediately
from Theorem \ref{T4M}.

The remainder of the paper after Section \ref{LOG} deals with what happens
when no normal convergence holds, i.e. when none of the equivalent
conditions in Theorem \ref{T3M} hold. Here, as is typical for studies of
non-normal convergence on\ Wiener space, the Malliavin calculus-based tools
such as Theorem \ref{T4M} are no longer useful. In such cases one may expect
that normalized centered quadratic variations of stationary Gaussian
sequences will converge in distribution towards second-chaos laws.

The classical result in this direction is a special case of what is often
known as the Dobrushin-Major--Taqqu theorem. This theorem refers to two
separate results from 1979 on non-linear functionals of fractional Brownian
motion (fBm). In the special case of centered quadratic variations, Taqqu 
(\cite{T79}) implies a convergence in the mean-square for the normalized
centered quadratic variations of fBm on $[0,1]$, towards the law of a
Rosenblatt r.v. $F_{\infty }$. On the other hand Dobrushin and Major (\cite{DM79}) prove
convergence to the same law, but only in law, not in the mean-square, for
the increments of fBm over intervals of length 1, as the horizon increases.
This second setting is the same as the one we use here: the increments of
fBm on unit intervals are known as fractional Gaussian noise (fGn), and as
mentioned above, fGn is the canonical power-scale example of a stationary
Gaussian sequence with unit variance. The self-similarity of fBm implies
that Taqqu's result and Dobrushin and Major's result are equivalent in the
setting we have chosen here, i.e. stationary Gaussian sequences. For the
historical reasons given above, it is legitimate to refer to such a result
as a Dobrushin-Major theorem.

Herein, we investigate the question of speed of convergence in the
Dobrushin-Major theorem. As we said, results of the type of Theorem \ref{T4M}
are of no use here, since the convergence we seek is of second-chaos type,
not normal. Instead, as in \cite{BN}, we refer to a result of Davydov and
Martinova \cite{DM}, which enables one to compare the total-variation
distance of two chaos variables with the square-root of their difference's $%
L^{2}$ norm. This requires the use of a change of probability space, which
is why the convergence cannot hold in $L^{2}\left( \Omega \right) $. In the
spirit of Dobrushin and Major, we consider the spectral representation of
the stationary Gaussian sequence, and obtain first a general strategy for
estimating the speed of convergence for sequences whose spectral density has
a certain type of functional asymptotic behavior under dilation of the
Fourier space's unit circle (Section \ref{SNN}). Then we specialize to cases
where the spectral density is a power times a slowly-varying function, in
order to present results similar to those for our log-modulated class of
examples. Our computational technique must go beyond the calculations in 
\cite{BN}, since the limit we obtain is the same for all log-modulated
processes for a fixed $H$, and therefore no self-similarity arguments may be
used. Our results in Section \ref{LOG2} show that the speed of convergence
itself does not depend on the modulation parameter: for modulations of the
form $\log ^{\beta }\left\vert x\right\vert $, $d_{TV}(F_{n},F_{\infty })=%
\mathcal{O}\left( \log ^{-2}n\right) $, for any $\beta >0$ (Theorem \ref%
{TNNS} on page \pageref{TNNS}).

We do not know if this result is sharp, but it does provide a sobering
extension to known speed-of-convergence results in the Dobrushin-Major
theorem, as compared for instance with the case $\beta =0$, where we get the
much faster power rate $d_{TV}(F_{n},F_{\infty })=\mathcal{O}\left(
n^{3/4-H}\right) $ (see Theorem 7.4.5 in \cite{npBook}). One can interpret
our result as saying that the price to pay for a mild \textquotedblleft
universality\textquotedblright\ result whereby a slowly-varying perturbation
of fGn still leads to a Rosenblatt limit in law, is that the speed of
convergence will be roughly as slow as the perturbation.

Summarizing the descriptions above, the remainder of this paper is
structured as follows. In Section \ref{CUM} we find sharp estimates of the
variance, third, and fourth cumulants of the normalized centered quadratic
variation $F_{n}$. We use these in Section \ref{THM} to prove a third-moment
theorem for normal convergence, with precise speed of convergence, for
general stationary Gaussian sequences (Theorem \ref{T3M}). In Section \ref%
{LOG}, we apply this theorem to a class of log-modulated covariance
structures, identifying a number of critical cases in the convergence rates
(Proposition \ref{logprop}), which go beyond the phase transitions recently
identified for fGn in \cite{alea}. In Section \ref{SNN}, we define a
strategy for estimating speeds of non-normal convergence for general $F_{n}$%
, based on a classical estimate in \cite{DM} and spectral representation of
stationary sequences (Corollary \ref{Cor2chaos}). In Section \ref{LOG2}, we
apply this strategy for sequences with power spectral density and log
modulation, proving that compared to the fGn case with $H>3/4$, the speed
deteriorates to $\log ^{-2}n$ instead of $n^{3/4-H}$ (Theorem \ref{TNNS}).

\section{Estimates for the second, third, and fourth cumulants\label{CUM}}

In this section we get precise estimates for the variance and the 3$^{rd}$
and 4$^{th}$ cumulants of $F_{n}$. Let us define first the comparison
relations that are used in what follows.

\begin{definition}[Comparison relations]
Given two deterministic numeric sequences $\left( a_{n}\right) _{n\geq
0},\left( b_{n}\right) _{n\geq 0}$ in a metric space, we use the following
notations and definitions for respectively domination, commensurability,
equivalence:

$a_{n}=\mathcal{O}\left( b_{n}\right) \Longleftrightarrow $ $\exists
C>0:a_{n}\leq Cb_{n}$ for $n$ large enough

$a_{n}\asymp b_{n}\Longleftrightarrow \exists c,C>0:cb_{n}\leq a_{n}\leq
Cb_{n}$ for $n$ large enough

$a_{n}\sim b_{n}\Longleftrightarrow \exists c_{n},C_{n}>0:\lim_{n\rightarrow
\infty }c_{n}=\lim_{n\rightarrow \infty }C_{n}=1$ and $c_{n}b_{n}\leq
a_{n}\leq C_{n}b_{n}$ for $n$ large enough.
\end{definition}

\begin{remark}
$\mathcal{O}$ is an order relation, while $\asymp $ and $\sim $ are
equivalence relations. Moreover, $a_{n}\mathcal{\asymp }b_{n}$ is equivalent
to \{$a_{n}=\mathcal{O}\left( b_{n}\right) $ and $b_{n}=\mathcal{O}\left(
a_{n}\right) $\}.
\end{remark}

\begin{remark}
\label{Cumul}The quantities $\kappa _{3}(F_{n}):=E(F_{n}^{3})$ and $\kappa
_{4}(F_{n}):=\mathbf{E}[F_{n}^{4}]-3$ are called the 3rd and 4th cumulants
of $F_{n}$. That $\kappa _{3}(F_{n})$ coincides with the third moment is
because $F_{n}$ is centered. Moreover, $\kappa _{4}(F_{n})$ is strictly
positive because $F_{n}$ is a non-Gaussian chaos r.v. (see \cite[Appendix A]%
{npBook} \ for details on cumulants on Wiener chaos).
\end{remark}

\begin{theorem}
\label{theo-est}Let $(F_{n})_{n\geq 0}$ be the sequence of normalized
centered quadratic variations of a centered stationary Gaussian sequence $%
\left( X_{n}\right) _{n\in \mathbf{Z}}$ with covariance $\rho $, i.e. with $%
\mathbf{E}\left[ X_{n}X_{n+k}\right] =\rho \left( k\right) $ for $n,k\in 
\mathbf{Z}$, we let $F_{n}:=V_{n}/\sqrt{v_{n}}$, where $V_{n}:=n^{-1/2}%
\sum_{k=0}^{n-1}(X_{k}^{2}-1)$, and $v_{n}:=E[V_{n}^{2}]$. Let $\kappa
_{3}(F_{n})$ and $\kappa _{4}(F_{n})$ be defined in Remark \ref{Cumul}.
Assume that the sequence of correlations $\rho $ has a constant sign, and
that $\left\vert \rho \right\vert $ is decreasing near $+\infty $, and $\rho
\left( 0\right) =1$. Then for large $n$,%
\begin{eqnarray}
\frac{1}{4}\frac{8}{v_{n}^{3/2}\sqrt{n}}\left( \sum_{|k|<n}\left\vert \rho
(k)\right\vert ^{3/2}\right) ^{2} &\leq &\left\vert \kappa
_{3}(F_{n})\right\vert \leq \frac{8}{v_{n}^{3/2}\sqrt{n}}\left(
\sum_{|k|<n}\left\vert \rho (k)\right\vert ^{3/2}\right) ^{2},
\label{equiv3} \\
\kappa _{4}(F_{n}) &=&\mathcal{O}\left( \frac{1}{v_{n}^{2}n}\left(
\sum_{|k|<n}\left\vert \rho (k)\right\vert ^{4/3}\right) ^{3}\right) ,
\label{equiv4} \\
v_{n} &=&-1+2\sum_{k=0}^{n-1}\left( 1-\frac{k}{n}\right) \rho ^{2}\left(
k\right) \asymp \sum_{k=0}^{n-1}\rho ^{2}\left( k\right) .  \label{equivv}
\end{eqnarray}
\end{theorem}

\noindent \emph{Proof.} We state the proof for the case of positive $\rho $;
for negative $\rho $, one only needs to replace $\rho $ by $\left\vert \rho
\right\vert $ in the computations.\emph{\vspace{0.1in}}

\noindent \emph{Step 1: Computation for the 3}$^{rd}$\emph{\ cumulant\vspace{%
0.1in}}

In \cite{4M} the following upper bound is proved: 
\begin{equation*}
\kappa _{3}(F_{n})\leq \frac{8}{v_{n}^{3/2}\sqrt{n}}\left( \sum_{|k|<n}\rho
(k)^{3/2}\right) ^{2}.
\end{equation*}%
That reference \cite{4M} contains the following explicit expression for the
third cumulant:%
\begin{equation*}
\kappa _{3}(F_{n})=\frac{8}{v_{n}^{3/2}\sqrt{n}}\left( \frac{1}{n}%
\sum_{j=0}^{n-1}\sum_{k,l=-j}^{n-1-j}\rho (k)\rho (k-l)\rho (l)\right) .
\end{equation*}%
By discarding the terms with $j>0$, we thus obtain the lower bound:%
\begin{equation*}
\kappa _{3}(F_{n})\geq \frac{8}{v_{n}^{3/2}\sqrt{n}}\left(
\sum_{k=0}^{n-1}\sum_{l=0}^{n-1}(1-\frac{k}{n})\rho (k)\rho (k-l)\rho
(l)\right) .
\end{equation*}%
Since $\rho $ is decreasing near $+\infty $, $\rho (k+l)\leq \sqrt{\rho
(k)\rho (l)}$ for large $k,l$. Thus we get 
\begin{eqnarray}
\kappa _{3}(F_{n}) &\geq &\frac{8}{v_{n}^{3/2}\sqrt{n}}\sum_{k=0}^{n-1}(1-%
\frac{k}{n})\rho (k)^{3/2}\sum_{l=0}^{n-1}\rho (l)^{3/2}  \notag \\
&\geq &\frac{8}{v_{n}^{3/2}\sqrt{n}}\sum_{l=0}^{n-1}\rho (l)^{3/2}\frac{1}{n}%
\sum_{k=0}^{n-1}S_{k}  \label{line2}
\end{eqnarray}%
where in line (\ref{line2}) we used an Abel summation given below in (\ref%
{abel}) (this relation is the subject of Step 2) on the sum $%
\sum_{k=0}^{n-1}(1-\frac{k}{n})\rho (k)^{3/2}$, and where $%
S_{n}=\sum_{l=0}^{n-1}\rho (l)^{3/2}$. However, since $\rho $ is positive
and decreasing, we have $S_{n}=\mathcal{O}\left( \frac{1}{n}%
\sum_{k=0}^{n-1}S_{k}\right) $ because $S_{2n}\leq 2S_{n}$ and $S_{2n}\geq
S_{n}$ for large $n$, so that $\frac{1}{n}\sum_{k=0}^{n-1}S_{n}\,\geq \frac{1%
}{2}S_{[\frac{n-1}{2}]}\,\geq \frac{1}{4}S_{n}$ for large $n$. That proves
the lower bound, and the Theorem's estimate on $\kappa _{3}$, modulo
estimate (\ref{abel}).\vspace{0.1in}

\noindent \emph{Step 2. Computation of the Abel summation.\vspace{0.1in}}

To compute 
\begin{equation*}
\sum_{k=0}^{n-1}(1-\frac{k}{n})\rho (k)^{3/2}=S_{n}-\frac{1}{n}%
\sum_{k=0}^{n-1}k\rho (k)^{3/2},
\end{equation*}%
we use Abel's summation-by-parts argument:%
\begin{eqnarray*}
\sum_{k=0}^{n-1}k\rho (k)^{3/2}
&=&\sum_{k=1}^{n-1}k(S_{k+1}-S_{k})=\sum_{k=1}^{n-1}(\left( k-1\right)
-k)S_{k}+(n-1)S_{n}+\rho (0)^{3/2} \\
&=&-\sum_{k=1}^{n-1}S_{k}+(n-1)S_{n}+\rho (0)^{3/2}
\end{eqnarray*}%
so that 
\begin{equation}
\sum_{k=0}^{n-1}(1-\frac{k}{n})\rho (k)^{3/2}=\frac{1}{n}(S_{n}-\rho
(0))^{3/2}+\sum_{k=1}^{n-1}S_{k}.  \label{abel}
\end{equation}%
\vspace{0.1in}

\noindent \emph{Step 3. Upper bound estimation for the 4}$^{th}$\emph{\
cumulant.\vspace{0.1in}}

In \cite{4M} the following upper bound is proved: $\kappa _{4}(F_{n})\leq 
\frac{c}{v_{n}^{2}n}\left( \sum_{|k|<n}\rho (k)^{4/3}\right) ^{3}.$ This is
our theorem's estimate on $\kappa _{4}\left( F_{n}\right) $. $\vspace*{0.1in}
$

\noindent \emph{Step 4. Asymptotic equivalent for the variance.\vspace*{0.1in%
}}

Let $T\left( n\right) :=\sum_{k=0}^{n-1}\rho ^{2}\left( k\right) $. Using a
change of summation variables, we find 
\begin{equation*}
nv_{n}=\sum_{k=0}^{n-1}\sum_{\ell =0}^{n-1}\rho \left( k-\ell \right)
^{2}=-n\rho ^{2}\left( 0\right) +2U\left( n\right)
\end{equation*}%
where 
\begin{equation*}
U\left( n\right) :=\sum_{k=0}^{n-1}\rho ^{2}\left( k\right) \left(
n-k\right) .
\end{equation*}%
Applying the same arguments as in steps 1 and 2, we find that $U\left(
n\right) \asymp nT\left( n\right) $. With $\rho \left( 0\right) =1$, this
finishes the proof of the theorem. $\hfill \blacksquare $ \bigskip

\section{A third-moment theorem for normal convergence\label{THM}}

In this section, we establish our third-moment theorem for the normalized
centered quadratic variation $F_{n}$ of our stationary Gaussian sequence $%
\left( X_{n}\right) $, with the assumptions as in Theorem \ref{theo-est}.

\begin{proposition}
\label{PropKey}Let $F_{n}$ and $\rho $ be as in Theorem \ref{theo-est}. Then 
$n^{1/4}\kappa _{4}(F_{n})^{3/4}=\mathcal{O}\left( \left\vert \kappa
_{3}(F_{n})\right\vert \right) $.
\end{proposition}

\noindent \emph{Proof :\vspace{0.1in}}

To lighten the notation, we assume that $\rho $ is positive. The estimates
obtained in Theorem \ref{theo-est} lead to 
\begin{equation*}
\frac{\kappa _{3}(F_{n})}{\kappa _{4}(F_{n})^{3/4}}\geq cn^{1/4}\left( \frac{%
\sum_{|k|<n}\rho (k)^{3/2}}{\left( \sum_{|k|<n}\rho (k)^{4/3}\right) ^{9/8}}%
\right) ^{2}.
\end{equation*}%
H\"{o}lder's inequality with $p=\frac{9}{4},q=\frac{9}{5}$ implies 
\begin{equation*}
\sum_{|k|<n}\rho (k)^{2/3}\rho (k)^{2/3}\leq \left( \sum_{|k|<n}\rho
(k)^{3/2}\right) ^{4/9}\left( \sum_{|k|<n}\rho (k)^{6/5}\right) ^{5/9}
\end{equation*}%
so that 
\begin{equation*}
\left( \frac{\sum_{|k|<n}\rho (k)^{3/2}}{\left( \sum_{|k|<n}\rho
(k)^{4/3}\right) ^{9/8}}\right) ^{2}\geq \left( \frac{\left(
\sum_{|k|<n}\rho (k)^{3/2}\right) ^{2/3}}{\left( \sum_{|k|<n}\rho
(k)^{6/5}\right) ^{5/6}}\right) ^{3/2}.
\end{equation*}%
Jensen inequality shows that $\left( \sum_{|k|<n}\rho (k)^{3/2}\right)
^{2/3}\geq \left( \sum_{|k|<n}\rho (k)^{6/5}\right) ^{5/6},$ which allows us
to conclude the proof. $\hfill \blacksquare $\bigskip

This proposition leads to the following.

\begin{theorem}[Third-moment theorem in total variation and convergence rates ]

\label{T3M}

Let $(F_{n})_{n\geq 0}$, $(v_{n})_{n\geq 0}$,and $\rho $ be as in Theorem %
\ref{theo-est}. Denote $N\sim \mathcal{N}(0,1)$. Then the following four
statements are equivalent:

\begin{description}
\item[(i)] $(F_{n})_{n\geq 0}$ converges in law towards $\mathcal{N}(0,1);$

\item[(ii)] $\lim_{n\rightarrow \infty }\mathbf{E}[F_{n}^{3}]=\mathbf{E}%
[N^{3}]=0;$

\item[(iii)] $\lim_{n\rightarrow \infty }\mathbf{E}[F_{n}^{4}]=\mathbf{E}%
[N^{4}]=3;$

\item[(iv)] $\left( \sum_{|k|<n}\left\vert \rho (k)\right\vert ^{3/2}\right)
^{2}=o\left( v_{n}^{3/2}\sqrt{n}\right) .$
\end{description}

\noindent Moreover the convergence rate in this situation is $\left\vert 
\mathbf{E}[F_{n}^{3}]\right\vert $, in the sense that for some $n_{0}>0$ 
\begin{equation*}
\exists c,C>0:\forall n>n_{0},\;c\left\vert \mathbf{E}[F_{n}^{3}]\right\vert
\leq d_{TV}(F_{n},N)\leq C\left\vert \mathbf{E}[F_{n}^{3}]\right\vert .
\end{equation*}%
We also have the following commensurability for this convergence rate:%
\begin{equation*}
\left\vert \mathbf{E}[F_{n}^{3}]\right\vert \asymp \frac{\left(
\sum_{|k|<n}\left\vert \rho (k)\right\vert ^{3/2}\right) ^{2}}{\left(
\sum_{|k|<n}\left\vert \rho (k)\right\vert ^{2}\right) ^{3/2}\sqrt{n}}.
\end{equation*}
\end{theorem}

\noindent \emph{Proof :\vspace{0.1in}}

By the 4th moment theorem (non-quantitative statement in Theorem \ref{T4M}),
(i) and (iii) are equivalent. By the commensurability relation (\ref{equiv3}%
) in Theorem \ref{theo-est}, (ii) and (iv) are equivalent. For the first
statement of the theorem, it is now enough to show (ii) and (iii) are
equivalent. By Proposition \ref{PropKey}, we have $\kappa _{4}(F_{n})=%
\mathcal{O}\left( \left\vert \kappa _{3}(F_{n})\right\vert
^{4/3}n^{-1/3}\right) $, thus if $\kappa _{3}(F_{n})\rightarrow 0$ then $%
\kappa _{4}(F_{n})\rightarrow 0$, i.e. (ii) implies (iii). On the other
hand, if (iii) holds, by Theorem \ref{T4M}, the convergence in (i) holds in
total variation, and the theorem's lower bound in (\ref{convergence_exacte})
implies (ii). The second statement of the theorem follows from (\ref%
{convergence_exacte}) and the estimate $\kappa _{4}\left( F_{n}\right) =%
\mathcal{O}\left( n^{-1/3}\left\vert \kappa _{3}\left( F_{n}\right)
\right\vert ^{4/3}\right) $ which shows that $\kappa _{4}\left( F_{n}\right)
=o\left( \left\vert \kappa _{3}\left( F_{n}\right) \right\vert \right) $.
The third statement combines relations (\ref{equiv3}) and (\ref{equivv}) in
Theorem \ref{theo-est}. The theorem is proved. $\hspace*{\fill}\blacksquare $%
\bigskip

The following result can be useful when dealing with non-normal convergence.

\begin{corollary}
\label{Cor-for-nonnormal}If $\rho \notin \ell ^{2}\left( \mathbf{Z}\right) $
then $v_{n}\sim 2\sum_{k=0}^{n-1}\left( 1-\frac{k}{n}\right) \rho ^{2}\left(
k\right) .$ This situation holds as soon as $(F_{n})_{n\geq 0}$ does not
converge in law towards $\mathcal{N}(0,1)$.
\end{corollary}

\noindent \emph{Proof :\vspace{0.1in}}

We use the notation and estimates from Step 4 in the proof of Theorem \ref%
{theo-est}. Since $U\left( n\right) \asymp nT\left( n\right) $, we have that 
$U\left( n\right) \asymp n$ iff $T\left( n\right) $ converges, i.e. $\rho
\in \ell ^{2}\left( \mathbf{Z}\right) $. Therefore if $\rho \notin \ell
^{2}\left( \mathbf{Z}\right) $ (i.e. $T\left( n\right) \rightarrow \infty $%
), then since $U\left( n\right) \geq \left( n/2\right) T\left( n/2\right)
\gg n$, we can ignore the first term in the equivalent (\ref{equivv}),
proving the asymptotic equivalent of the corollary.

To prove the corollary's second statement, by Theorem \ref{T3M}, it is
sufficient to show that $\rho \in \ell ^{2}\left( \mathbf{Z}\right) $
implies condition (iv) of the theorem. Since $\rho \in \ell ^{2}\left( 
\mathbf{Z}\right) $ implies that $v_{n}$ is bounded, we only need to show
that $S\left( n\right) :=\sum_{|k|<n}\left\vert \rho (k)\right\vert
^{3/2}=o\left( n^{1/4}\right) $. The case of $S\left( n\right) $ bounded is
trivial, so we assume $S\left( n\right) $ is unbounded. For any fixed $%
\varepsilon >0$, since $\rho \in \ell ^{2}$, there exists $n_{1}$ such that $%
\sum_{k=n_{1}}^{\infty }\rho \left( k\right) ^{2}\leq \varepsilon $. Also,
since $S\left( n_{1}\right) $ is fixed, and $S\left( n\right) $ diverges,
there exists $n_{2}>n_{1}$ such that for all $n\geq n_{2}$, $S\left(
n\right) \geq 2S\left( n_{1}\right) $. Hence we have 
\begin{equation*}
S\left( n\right) =S\left( n\right) -S\left( n_{1}\right) +S\left(
n_{1}\right) \leq S\left( n\right) -S\left( n_{1}\right) +2^{-1}S\left(
n\right) 
\end{equation*}%
so that%
\begin{equation*}
S\left( n\right) \leq 2\left( S\left( n\right) -S\left( n_{1}\right) \right)
.
\end{equation*}%
We have by Jensen's inequality, for all $n\geq n_{2}$,%
\begin{eqnarray*}
S\left( n\right) -S\left( n_{1}\right)  &=&2\sum_{k=n_{1}}^{n-1}\left\vert
\rho (k)\right\vert ^{3/2}=2\left( n-n_{1}\right) \left[ \left(
\sum_{k=n_{1}}^{n-1}\frac{1}{n-n_{1}}\left\vert \rho (k)\right\vert
^{3/2}\right) ^{4/3}\right] ^{3/4} \\
&\leq &2\left( n-n_{1}\right) ^{1/4}\left[ \sum_{k=n_{1}}^{n-1}\frac{1}{%
n-n_{1}}\left\vert \rho (k)\right\vert ^{2}\right] ^{3/4}\leq 2\varepsilon
^{3/4}\left( n-n_{1}\right) ^{1/4}.
\end{eqnarray*}%
This finishes the proof of the corollary. $\hfill \blacksquare $\bigskip 

\section{Example: sequences with log--modulated power covariance\label{LOG}}

Assume here that the correlation function $\rho $ satisfies%
\begin{equation}
\left\vert \rho (n)\right\vert \sim n^{2H-2}\log ^{2\beta }(n),
\label{rholog}
\end{equation}%
for large $n$, with $H\in \lbrack 0,1]$ and $\beta \in \mathbf{R}$. A
stationary sequence $\left( X_{n}\right) _{n\in \mathbf{Z}}$ with this
property can easily be constructed as a Gaussian Fourier integral 
\begin{equation*}
X_{n}=\int_{S_{1}}\sqrt{q\left( x\right)} \cos \left( xn\right) W\left(
dx\right) +\int_{S_{1}}\sqrt{q\left( x\right)} \sin \left( xn\right) \tilde{W%
}\left( dx\right) 
\end{equation*}%
with Fourier coefficients $q\left( x\right) =x^{1-2H}\log ^{-2\beta }(x^{-1})
$, where $x$ is the Fourier parameter on the unit circle $S^{1}$, and $W$
and $\tilde{W}$ are independent white noises on $S^{1}$. The details are
omitted here; some can be found in the Appendix (see Section \ref{Complex}).
We also assume, as in the past, that $\rho $ is of constant sign and $%
\left\vert \rho \right\vert $ decreases for large $n$. The process
constructed above by Gaussian Fourier integral satisfies this assumption. A
arbitrary constant scaling factor should be added to the asymptotics of $%
\left\vert \rho \right\vert $ (like the $H\left\vert 2H-1\right\vert $ for
fBm), but we omit this in this section for notational simplicity.

The case of discrete-time fractional Gaussian noise (fGn), treated in \cite%
{4M} and \cite{alea}, falls within the special case of $H\in \left(
0,3/4\right) $ and $\beta =0$. An extension of the fGn is developed in \cite%
{BBL} contains a correction term added to the the power spectral density of
fractional Brownian motion, which decays like a faster power at high
frequency. This means that the additional term is treated like a lower-order
remainder, for continuous-time data. Strictly speaking, this type of model
is not immediately comparable to the discrete-time models we consider here.
However, a similar study for discrete time could be done, for instance, on
correlation structures of the form $\left\vert \rho (n)\right\vert \sim
n^{2H-2}+r\left( n\right) $ where the remainder $r\left( n\right) =O\left(
n^{2H-2-\gamma }\right) $ for some $\gamma >0$. Such a class is contained in
our assumption (\ref{rholog}) with $\beta =0$; the corresponding convergence
results and speeds of convergence would then depend only on $H$, not $\gamma 
$. We omit any further discussion of this point for the sake of brevity.

We will see below in Proposition \ref{logprop}, that we can cover the case $%
H=3/4$ for any $\beta $, and that within the two presumed \textquotedblleft
critical thresholds\textquotedblright\ $H=2/3$ and $H=3/4$, there arise
further \textquotedblleft critical log-thresholds\textquotedblright\ $\beta
=-1/3$ and $\beta =-1/4$ respectively; moreover, the thresholds $H=2/3$ and $%
H=3/4$ only give rise to \textquotedblleft exotic\textquotedblright\
convergence rates (i.e. with \textquotedblleft log
corrections\textquotedblright ) for certain ranges of the parameter $\beta $%
. In this sense, the notion of critical value or of phase transition is
model-class-dependent. The reader could further convince herself of this by
considering a class of processes with $H=2/3$, $\beta =-1/3$, and an
additional factor $\log \log ^{2\gamma }\left( n\right) $, to find out that
the \textquotedblleft critical pair\textquotedblright\ $\left( H=2/3,\beta
=-1/3\right) $ which we exhibit harbors further ranges and cutoff values of
the parameter $\gamma $, some of which may be considered more exotic than
others. In other words, the so-called critical cases do not have any
fundamental significance, but just appear as consequences of the models'
scaling choices (power-scale in the fractional Brownian example,
log+power-scale in our class of examples).

Also notice that by Proposition \ref{PropKey}, the asymptotics of $%
d_{TV}(F_{n},N)$ are always given by those of $\kappa _{3}\left(
F_{n}\right) $, even when the Breuer-Major theorem fails. This is not a
robust result, however. For instance, one can check from the calculations in 
\cite{alea} and \cite{4M} that for $q\geq 3$ and certain values of $H$, the
third and fourth cumulants take turns at determining the speed of
convergence of $q$th-power variations of fBm. As soon as the limit of $F_{n}$
is not normal, the question of determining an optimal speed of convergence
for $F_{n}$ becomes unresolved. It is known that for fBm and other
self-similar processes, when $H>3/4$, the normalized quadratic variation
converges to a so-called Rosenblatt distribution, which is a law in the
second chaos, which depends on the parameter $H$ (see \cite{fvRosen} and
references therein). Estimating the rate of this convergence from above, for
fBm and for log-modulated processes, is the topic of Sections \ref{SNN} and %
\ref{LOG2}; therein we will see that unlike the case of normal convergence,
the rate is determined by the modulation rather than $H$.

Returning to the topic of normal convergence of $F_{n}$, under the
asymptotics in (\ref{rholog}) we now compute the equivalents of the
convergence rates exactly, thanks to Theorem \ref{theo-est}. Let us recall
the following well-known result about Bertrand series.

\begin{pty}[Equivalents of Bertrand series]
The series $S_{n}(\alpha ,\beta ):=\sum_{n>0}n^{\alpha }\log ^{\beta }(n)$
converges if and only if $\alpha <-1$ or $\alpha =-1$ and $\beta <-1$. When
the series diverges, we have the following equivalents for its partial sum:

\begin{itemize}
\item $S_{n}(-1,-1)\sim \log (\log (n));$

\item $S_{n}(-1,\beta )\sim \frac{1}{\beta +1}\log ^{\beta +1}(n)$ if $\beta
>-1;$

\item $S_{n}(\alpha ,\beta )\sim \frac{1}{\alpha +1}n^{\alpha +1}\log
^{\beta }(n)$ if $\alpha >-1,\beta >-1.$
\end{itemize}
\end{pty}

We may now use this lemma and Theorem \ref{theo-est} to obtain the
asymptotic order of the third and fourth cumulants, which gives us the rate
of convergence to the normal law for $\left( F_{n}\right) $, as a
consequence of Theorem \ref{T3M}.

\begin{proposition}
\label{logprop}With the notation and assumptions as in Theorem \ref{T3M},
there are positive constants $c,C$ such that $cM_{n}\leq d_{TV}(F_{n},N)\leq
CM_{n}$ where

\begin{itemize}
\item $M_{n}=\frac{1}{\sqrt{n}}$ if $H<\frac{2}{3}$ or $H=\frac{2}{3},\beta
<-\frac{1}{3},$

\item $M_{n}=\frac{\log (\log (n))^{2}}{\sqrt{n}}$ if $H=\frac{2}{3},\beta =-%
\frac{1}{3},$

\item $M_{n}=\frac{1}{\sqrt{n}}\log ^{2(3\beta +1)}(n)$ if $H=\frac{2}{3}%
,\beta >-\frac{1}{3},$

\item $M_{n}=n^{6H-\frac{9}{2}}\log ^{6\beta }(n)$ if $\frac{2}{3}<H<\frac{3%
}{4}$ or $H=\frac{3}{4},\beta <-\frac{1}{4},$

\item $M_{n}=\log ^{-\frac{3}{2}}(n)\log (\log (n))^{-\frac{3}{2}}$ if $H=%
\frac{3}{4},\beta =-\frac{1}{4},$

\item $M_{n}=\log ^{-\frac{3}{2}}(n)$ if $H=\frac{3}{4},\beta >-\frac{1}{4}.$
\end{itemize}

Since $M_{n}\rightarrow 0$ in all these cases, $d_{TV}(F_{n},N)\rightarrow 0$
at the same rates. The normalizing factor $v_{n}$ converges if and only if \{%
$H<\frac{3}{4}$ or $H=\frac{3}{4},\beta <-\frac{1}{4}$\}.
\end{proposition}

\begin{remark}
If \{$H=\frac{3}{4},\beta \geq -\frac{1}{4}$\} or if $H>3/4$, the
Breuer-Major theorem fails by definition because $V_{n}:=n^{-1/2}%
\sum_{k=0}^{n-1}(X_{k}^{2}-1)$ does not converge to a normal (its variance
diverges). However the last two cases in the previous proposition show that
normal convergence still holds for all cases where $H=3/4$. In other words,
one must distinguish between (i) a Breuer-Major-type theorem, which attempts
to characterize situations where a central limit theorem might hold for
partial sums of highly dependent sequences $X_{k}^{2}-1$, where the familiar
normalization $n^{-1/2}$ can still be used, and (ii) other normal
convergences where the dependence of the terms $X_{k}^{2}-1$ is too strong
for a central-limit normalization, but hold under a stronger normalization.
\end{remark}

\noindent \emph{Proof of Proposition \ref{logprop}.\vspace{0.1in}}

All the computations below are based on the sharp estimates of Theorem \ref%
{theo-est}. By Theorem \ref{T3M}, it is sufficient to compute a
commensurable equivalent of $\kappa _{3}(F_{n})$.\emph{\vspace{0.1in}}

\noindent \emph{Step 1: computing the series }$\sum_{|k|<n}\rho (k)^{3/2}$.%
\vspace{0.1in}

Using our class of examples for $\rho $, we have that $\sum_{|k|<n}\rho
(k)^{3/2}$ diverges if and only if $\sum_{|k|<n}|k|^{3(H-1)}\log ^{3\beta
}(|k|)$ diverges, and is equivalent to this divergent partial sum in that
case. Therefore, by the equivalents for Bertrand series, $\sum_{|k|<n}\rho
(k)^{3/2}$ is equivalent to

\begin{description}
\item[$\ast $] a constant $l(H,\beta )>0$ if $H<\frac{2}{3}$ or $H=\frac{2}{3%
},\beta <-\frac{1}{3}$

\item[$\ast $] $\log (\log (n))$ if $H=\frac{2}{3},\beta =-\frac{1}{3}$

\item[$\ast $] $\frac{1}{3\beta +1}\log ^{3\beta +1}n$ if $H=\frac{2}{3}%
,\beta >-\frac{1}{3}$

\item[$\ast $] $\frac{1}{3H-2}n^{3H-2}\log ^{3\beta }n$ if $H>\frac{2}{3}$
\end{description}

\noindent \emph{Step 2: computing }$v_{n}.$\vspace{0.1in}

Again, by definition of $\rho $, we have that the following quantities are
divergent simultaneously and the following equivalences hold in that case: $%
v_{n}\asymp \sum_{|k|<n}\rho (k)^{2}\sim \sum_{|k|<n}|k|^{4(H-1)}\log
^{4\beta }(|k|)$. Therefore $v_{n}$ is equivalent to

\begin{description}
\item[$\ast $] a constant $l^{\prime }(H,\beta )>0$\ if $H<\frac{3}{4}$ or $%
H=\frac{3}{4},\beta <-\frac{1}{4}$

\item[$\ast $] $\log (\log (n))$ if $H=\frac{3}{4},\beta =-\frac{1}{4}$

\item[$\ast $] $\frac{1}{3\beta +1}\log ^{4\beta +1}n$ if $H=\frac{3}{4}%
,\beta >-\frac{1}{4}$

\item[$\ast $] $\frac{1}{4H-3}n^{4H-3}\log ^{4\beta }n$ if $H>\frac{3}{4}$
\end{description}

\noindent \emph{Step 3: computing }$\kappa _{3}(F_{n})$.\vspace*{0.1in}

We can now compute a commensurable equivalent for $\kappa _{3}(F_{n})$
thanks to relation (\ref{equiv3}).

\begin{description}
\item[$\ast $] $\kappa _{3}(F_{n})\asymp \frac{1}{\sqrt{n}}$ if $H<\frac{2}{3%
}$ or $H=\frac{2}{3},\beta <-\frac{1}{3}$

\item[$\ast $] $\kappa _{3}(F_{n})\asymp \frac{\log (\log (n))^{2}}{\sqrt{n}}
$ if $H=\frac{2}{3},\beta =-\frac{1}{3}$

\item[$\ast $] $\kappa _{3}(F_{n})\asymp \frac{1}{\sqrt{n}}\log ^{2(3\beta
+1)}(n)$ if $H=\frac{2}{3},\beta >-\frac{1}{3}$

\item[$\ast $] $\kappa _{3}(F_{n})\asymp n^{6H-\frac{9}{2}}\log ^{6\beta
}(n) $ if $\frac{2}{3}<H<\frac{3}{4}$ or $H=\frac{3}{4},\beta <-\frac{1}{4}$

\item[$\ast $] $\kappa _{3}(F_{n})\asymp \log ^{-\frac{3}{2}}(n)\log (\log
(n))^{-\frac{3}{2}}$ if $H=\frac{3}{4},\beta =-\frac{1}{4}$

\item[$\ast $] $\kappa _{3}(F_{n})\asymp \log ^{-\frac{3}{2}}(n)$ if $H=%
\frac{3}{4},\beta >-\frac{1}{4}$
\end{description}

Theorem \ref{T3M} now allows us to conclude.\hfill $\blacksquare $

\section{Strategy for non-normal convergence\label{SNN}}

In Theorem \ref{T3M}, we saw that $\left( F_{n}\right) _{n\geq 0}$ converges
to a normal if and only if Condition (iv) therein is satisfied. When this
condition does not hold, one may wonder what kind of other convergence we
could get. The celebrated theorem of Dobrushin-Major and Taqqu \cite{DM79,
T79} spells out the possible second-chaos limits of $\left( F_{n}\right)
_{n\geq 0}$ which occur for slowly-varying perturbations of fGn. In this
section, we detail a general methodology, based on the context used by
Dobrushin and Major, to determine whether convergence to a second-chaos
limit holds.

We are not able to provide general criteria as sharp and as explicit as
those which we gave in the previous section for normal convergence; we are
not aware of any such works in the literature. However, a classical result
of Davydov and Martinova \cite{DM} enables us to define a strategy for
establishing speed of convergence in total variation to second-chaos limits.
Section \ref{LOG2} shows how to apply this methodology when the
slowly-varying terms are logarithmic, similarly to the class of examples
presented in Section \ref{LOG}.

It was recently established in \cite{NPoly} that the only limits in
distribution for sequences in a 2nd Wiener chaos are of the form $\mathcal{L}%
\left( N+F\right) $ where $F$ and $N$ are independent and $F$ is the law of
a second-chaos rv, and $N$ is Gaussian. It is conceivable that one could
find a choice of law for our sequence $X$ such that $\left( F_{n}\right)
_{n\geq 0}$ converges in law to such a convolution. We will avoid such a
situation, and investigate how to prove instead that $\left( F_{n}\right)
_{n\geq 0}$ converges in law to the law of a second-chaos random variable.

As in \cite{T79}, we can express our stationary sequence $X$ as a Fourier
integral over the unit circle $S^{1}\equiv \lbrack -\pi ,\pi )$. Let $q$ be
the Fourier transform of the even sequence $\rho $ on $\mathbf{Z}$, defined
on $S^{1}$ by 
\begin{equation*}
q\left( x\right) :=\sum_{n\in \mathbf{Z}}\rho \left( n\right) \cos \left(
nx\right) .
\end{equation*}%
Assume that $q\in L^{1}\left( S^{1}\right) $. It is known that $q$ is
non-negative, with $\int_{S^{1}}q\left( x\right) dx=\rho \left( 0\right) =1$%
. Since we are interested in sequences $X$ such that $\left( F_{n}\right)
_{n\geq 0}$ does not converge to a normal, we will find that typically we
have a long memory property, i.e. $\rho \notin \ell ^{1}\left( \mathbf{Z}%
\right) $, so that the Fourier series defining $q$ is not absolutely
convergent. In fact, for $\rho $ decreasing, we saw in Corollary \ref%
{Cor-for-nonnormal} that $\rho $ is not even in $\ell ^{2}\left( \mathbf{Z}%
\right) $. However, if $q\in L^{1}\left( S^{1}\right) $ is given
exogenously, the classical theorem on pointwise convergence of Fourier
series (see \cite{FourierBook} Chapter 3, Theorem 2.1) implies that the
Fourier series of $q$ converges to $q$ at all points where $q$ is
differentiable. The assumption $q\in L^{1}\left( S^{1}\right) $ and its
differentiability will be verified in our examples in Section \ref{LOG2}. We
will also need to ensure that $\rho $ satisfies a Fourier inversion theorem.
Since $\rho \notin \ell ^{1}\left( \mathbf{Z}\right) $, we cannot appeal
directly to the classical Fourier inversion theorem, and we must thus check
that it holds on a case-by-case basis. In Section \ref{LOG2}, we will work
with covariances $\rho $ which are defined in terms of their Fourier
transform $q$, so that Fourier inversion is automatic.

With $q\in L^{1}\left( S^{1}\right) $ as above, and assuming that Fourier
inversion holds for $\rho $, then as for instance in Dobrushin \cite%
{DubroRGF}, the following spectral representation holds for the centered
Gaussian sequence $X$ whose covariance is $\rho $: there exists a standard
complex-valued white noise $W$ on $S^{1}$ such that%
\begin{equation}
X(k)=\int_{S^{1}}e^{ikx}\sqrt{q(x)}W(dx).  \label{SR}
\end{equation}%
This appellation for $W$ means that for $x\in \lbrack 0,\pi ]$, $W\left(
dx\right) =B_{1}\left( dx\right) +iB_{2}\left( dx\right) $, where $B_{1}$
and $B_{2}$ are two real-valued independent white noise measures on $[0,\pi
] $ with scaling constant $(2\pi )^{-1/2}$ (i.e. $Var[B_{i}([0,\pi ])]=1/2$%
), and for every $x\in \lbrack 0,\pi ]$, $W\left( -dx\right) =\overline{%
W\left( dx\right) }=B_{1}\left( dx\right) -iB_{2}\left( dx\right) $. It is
helpful to realize that the representation (\ref{SR}) is equivalent to 
\begin{equation*}
X\left( k\right) =\int_{-\pi }^{\pi }\sqrt{q\left( x\right) }\cos \left(
kx\right) W_{1}\left( dx\right) +\int_{-\pi }^{\pi }\sqrt{q\left( x\right) }%
\sin \left( kx\right) W_{2}\left( dx\right)
\end{equation*}%
where $W_{1}$ and $W_{2}$ are i.i.d. real-valued white noises on $[-\pi ,\pi
]$, standardized so that $Var\left( W_{i}\left[ -\pi ,\pi \right] \right) =1$%
. Details on the properties of $W$ are given in the Appendix in Section \ref%
{Complex}. In the remainder of this section and in Section \ref{LOG2}, it
will be sufficient for us to know the following special case of the isometry
formula for double Wiener integrals with respect to $W$: for $f\in
L^{2}\left( \left( S^{1}\right) ^{2},\mathbf{C}\right) $, if $f$ satisfies
the Hermitian evenness property $f\left( -x,-y\right) =\overline{f\left(
x,y\right) }$, then%
\begin{equation}
\mathbf{E}\left[ I_{2}\left( f\right) ^{2}\right] =\iint_{[-\pi ,\pi
]^{2}}\left\vert f\left( x,y\right) \right\vert ^{2}\frac{dxdy}{\left( 2\pi
\right) ^{2}}=:\left\Vert f\right\Vert _{L^{2}\left( \left( S^{1}\right)
^{2}\right) }^{2}.  \label{isometry}
\end{equation}%
See Section \ref{Complex} for a proof. Note in addition that the isometry
formula (\ref{isometry}) for double Wiener integrals against a
complex-valued white noise does not contain the usual factor of $2$
associated with the isometry property for double Wiener integrals with
respect to a real-valued white noise.

From formula (\ref{SR}) and the product formula for Wiener integrals, we can
write 
\begin{equation*}
X\left( k\right) ^{2}-1=\iint_{\left( S^{1}\right) ^{2}}W\left( dx\right)
W\left( dy\right) e^{ik\left( x+y\right) }\sqrt{q\left( x\right) q\left(
y\right) }.
\end{equation*}%
Using the stochastic Fubini theorem, justified because $q\in L^{1}\left(
S^{1}\right) $ implies that $\left( k,x,y\right) \mapsto e^{ik\left(
x+y\right) }\sqrt{q\left( x\right) q\left( y\right) }$ is in $L^{2}\left(
\left( S^{1}\right) ^{2}\right) \times \ell ^{2}\left( \mathbf{Z}\right) $,
we get%
\begin{equation*}
F_{n}=\frac{1}{\sqrt{nv_{n}}}\iint_{\left( S^{1}\right) ^{2}}W\left(
dx\right) W\left( dy\right) \frac{e^{in\left( x+y\right) }-1}{e^{i\left(
x+y\right) }-1}\sqrt{q\left( x\right) q\left( y\right) }.
\end{equation*}

To prove that $\left( F_{n}\right) _{n\geq 0}$ converges in law to a
second-chaos distribution, it is sufficient to prove that there exists a
Wiener space $\left( \hat{\Omega},\mathcal{\hat{F}},\mathbf{\hat{P}}\right) $
such that for every $n$ there is a random variable $\hat{F}_{n}$ with the
same law as $F_{n}$, and another second-chaos random variable $F_{\infty }$
on $\left( \hat{\Omega},\mathcal{\hat{F}},\mathbf{\hat{P}}\right) $ (not
dependent on $n$), such that $\mathbf{\hat{E}}\left[ \left( F_{\infty }-\hat{%
F}_{n}\right) ^{2}\right] \rightarrow 0$. Furthermore, since $d\left(
F_{\infty },F_{n}\right) =d\left( F_{\infty },\hat{F}_{n}\right) $ for any
distance $d$ on the set of laws, to estimate the total-variation distance
between the law of $F_{\infty }$ and the law of $F_{n}$, one may rely on a
theorem of Davydov and Martinova \cite{DM}, by which, if indeed $\left( \hat{%
F}_{n}\right) _{n>0}$ converges in $L^{2}\left( \Omega \right) $ to $%
F_{\infty }$, then%
\begin{equation}
d_{TV}\left( \hat{F}_{n},F_{\infty }\right) \leq c_{F_{\infty }}\left( 
\mathbf{\hat{E}}\left[ \left( \hat{F}_{n}-F_{\infty }\right) ^{2}\right]
\right) ^{1/4}  \label{DMineq}
\end{equation}%
where $c_{F_{\infty }}$ is a finite constant depending only on the law of $%
F_{\infty }$.

We now change variables from $\left( x,y\right) $ to $\left( x^{\prime
},y^{\prime }\right) :=\left( xn,yn\right) $, omitting the primes for
parsimony of notation, and we write $I_{n}:=[-\pi n,\pi n)$ for the
corresponding scaled circle. The self-similarity of $W\left( dx\right) $
with index $1/2$ means that $W\left( n^{-1}dx\right) $ has the same law as $%
n^{-1/2}W\left( dx\right) $. Therefore the random variable $F_{n}$ has the
same law as the variable $\hat{F}_{n}$ defined as follows under a standard
complex-valued white noise measure $\hat{W}$ on $\mathbf{R}$ scaled by $%
\left( 2\pi \right) ^{-1/2}$ (in particular $\left( \hat{\Omega},\mathcal{%
\hat{F}},\mathbf{\hat{P}}\right) $ is the $\mathbf{C}$-valued Wiener space
of $\hat{W}$, and $\left\vert \hat{W}\left( dx\right) \right\vert
^{2}=dx/\left( 2\pi \right) $ ): 
\begin{equation*}
\hat{F}_{n}=\frac{1}{\sqrt{nv_{n}}}\iint_{\mathbf{R}^{2}}\mathbf{1}%
_{I_{n}}\left( x\right) \mathbf{1}_{I_{n}}\left( y\right) \sqrt{q\left( 
\frac{x}{n}\right) q\left( \frac{y}{n}\right) }\frac{1/n}{e^{i\left(
x+y\right) /n}-1}\left( e^{i\left( x+y\right) }-1\right) \hat{W}\left(
dx\right) \hat{W}\left( dy\right) .
\end{equation*}

We note that $g_{n}:x\mapsto \left( 1/n\right) \left( e^{ix}-1\right)
/\left( e^{inx}-1\right) $ converges pointwise and boundedly to the bounded
function $g:x\mapsto \left( e^{ix}-1\right) /\left( ix\right) $, including
at $x=0$; this is proved using elementary calculations and the inequality $%
\left\vert e^{ix}-1-ix\right\vert \leq \left\vert x\right\vert ^{2}$ for $%
\left\vert x\right\vert \leq 1/2$. This implies by dominated (bounded)
convergence that $g_{n}$ and $g$ can be interchanged in expressions which
are $L^{2}$-convergent. We also recall that we may use the expression $%
nv_{n}=-n+2\sum_{k=0}^{n-1}\left( n-k\right) \rho ^{2}\left( k\right) $
which was established in (\ref{equivv}) in Theorem \ref{theo-est}.

Now using estimate (\ref{DMineq}) and the isometry property (\ref{isometry}%
), the above discussion proves the following criterion for convergence in
law of $\hat{F}_{n}$, which includes a possible quantitative estimate. We
emphasize that the convergence-in-law portion of this strategy is not new,
since it was used by Dobrushin and Major in the case of fGn with slowly
varying modulation.

\begin{theo}[Rate of convergence {in total variation to the Rosenblatt law}]

\label{T2chaos}Let $(F_{n})_{n\geq 0}$ be the sequence of normalized
quadratic variations of a centered stationary Gaussian sequence $\left(
X_{n}\right) $ with covariance $\rho $, i.e. $F_{n}:=V_{n}/\sqrt{v_{n}}$,
where $V_{n}:=n^{-1/2}\sum_{k=0}^{n-1}(X_{k}^{2}-1)$, and $%
v_{n}:=E[V_{n}^{2}]=n^{-1}\sum_{k=0}^{n-1}\sum_{\ell =0}^{n-1}\rho \left(
k-\ell \right) ^{2}=-1+2\sum_{k=0}^{n-1}\left( 1-k/n\right) \rho ^{2}\left(
k\right) $. Assume that the spectral density $q$ of $X$, defined by%
\begin{equation*}
q\left( x\right) :=\sum_{n\in \mathbf{Z}}\rho \left( n\right) \cos \left(
nx\right) 
\end{equation*}%
exists as a member of $L^{1}\left( [-\pi ,\pi )\right) $ Assume that Fourier
inversion holds, i.e. that $\int_{-\pi }^{\pi }e^{ikx}q\left( x\right) \frac{%
dx}{2\pi }=\rho \left( k\right) $ for all $k\in \mathbf{Z}$. Assume $\rho
\left( 0\right) =1$. Let 
\begin{equation*}
I_{n}:=[-\pi n,\pi n),\;\;g\left( x\right) :=\frac{e^{ix}-1}{ix}.
\end{equation*}%
Assume there exists a real function $f\in L_{loc}^{2}\left( \mathbf{R}%
^{2}\right) $ which is even in both variables and such that $\left(
x,y\right) \mapsto \left( nv_{n}\right) ^{-1/2}\sqrt{q\left( \frac{x}{n}%
\right) q\left( \frac{y}{n}\right) }\mathbf{1}_{I_{n}}\left( x\right) 
\mathbf{1}_{I_{n}}\left( y\right) g\left( x+y\right) $ converges in $%
L^{2}\left( \mathbf{R}^{2}\right) $ to $f\left( x,y\right) g\left(
x+y\right) $.

Then $F_{n}$ converges in law to the law of a second-chaos variable $%
F_{\infty }$, and we have the representation%
\begin{equation*}
F_{\infty }=\iint_{\mathbf{R}^{2}}f\left( x,y\right) \frac{e^{i(x+y)}-1}{%
i\left( x+y\right) }W\left( dx\right) W\left( dy\right)
\end{equation*}%
where $W$ is a complex white noise on $\mathbf{R}$ with scale determined by $%
\left\vert W\left( dx\right) \right\vert ^{2}=dx/\left( 2\pi \right) $.

Moreover, with $c_{F_{\infty }}$ the constant in (\ref{DMineq}), the speed
of convergence in total variation is bounded above as%
\begin{eqnarray*}
&&d_{TV}\left( F_{n},F_{\infty }\right) ^{4} \\
&\leq &\left( c_{F_{\infty }}\right) ^{4}\iint_{\mathbf{R}^{2}}\left\vert 
\mathbf{1}_{I_{n}}\left( x\right) \mathbf{1}_{I_{n}}\left( y\right) \frac{1}{%
\sqrt{nv_{n}}}\sqrt{q\left( \frac{x}{n}\right) q\left( \frac{y}{n}\right) }%
\frac{n^{-1}\left( e^{i\left( x+y\right) }-1\right) }{e^{i\left( x+y\right)
/n}-1}-f\left( x,y\right) \frac{e^{i(x+y)}-1}{i\left( x+y\right) }%
\right\vert ^{2}dxdy
\end{eqnarray*}
\end{theo}

\noindent \emph{Proof.} See above development. \hfill $\blacksquare $

\begin{rem}
The examples in Section \ref{LOG2} satisfy the assumption on $q$ and $\rho $
in Theorem \ref{T2chaos}.
\end{rem}

\begin{rem}
\label{rem-even}The last inequality above, and Corollary \ref{Cor2chaos}
below, hold with $y$ replaced by $-y$ since $f$ is even in each variable,
and $q$ and $I_{n}$ are even.
\end{rem}

\begin{rem}
Under the assumptions of Theorem \ref{T2chaos}, the law of $\left(
F_{n}\right) _{n>0}$ cannot converge to a normal law. Consequently by
Theorem \ref{T3M}, the other three equivalent conditions (ii), (iii), (iv)
therein fail. In particular, by Corollary \ref{Cor-for-nonnormal}, $\rho
\notin \ell ^{2}\left( \mathbf{Z}\right) $.
\end{rem}

The following corollary is useful to estimate the speed of convergence in
Theorem \ref{T2chaos}. It enables one to introduce a trade-off between the
speed of convergence of the improper integral defining $\left\Vert
fg\right\Vert _{L^{2}\left( \mathbf{R}^{2}\right) }^{2}$ and the speed of
convergence of $\sqrt{q\left( \frac{x}{n}\right) q\left( \frac{y}{n}\right)
n^{-1}/v_{n}}$ to $f$.

\begin{cor}
\label{Cor2chaos}Under the assumptions and notation of Theorem \ref{T2chaos}%
, for any $\alpha \in \left( 0,1\right) $, with $h\left( x\right) :=\min
\left( 1,\left\vert x\right\vert ^{-1}\right) $, and $nv_{n}=-n+2%
\sum_{k=0}^{n-1}\left( n-k\right) \rho ^{2}\left( k\right) $,%
\begin{eqnarray}
\frac{1}{\left( c_{F_{\infty }}\right) ^{4}}d_{TV}\left( F_{n},F_{\infty
}\right) ^{4} &\leq &5\iint_{\mathbf{R}^{2}\backslash \left( I_{n^{\alpha
}}\right) ^{2}}\left\vert f\left( x,y\right) \right\vert ^{2}h^{2}\left(
x+y\right) dxdy  \label{cor2chaosline1} \\
&&+16\iint_{\left( I_{n}\right) ^{2}}\left\vert \sqrt{\frac{q\left( \frac{x}{%
n}\right) q\left( \frac{y}{n}\right) }{nv_{n}}}-f\left( x,y\right)
\right\vert ^{2}h^{2}\left( x+y\right) dxdy  \label{cor2chaosline2} \\
&&+\frac{4}{n^{2-2\alpha }}\iint_{\mathbf{R}^{2}}\left\vert f\left(
x,y\right) \right\vert ^{2}h^{2}\left( x+y\right) dxdy
\label{cor2chaosline3}
\end{eqnarray}
\end{cor}

\noindent \emph{Proof.} Write $g\left( x\right) =\left( e^{ix}-1\right)
/\left( ix\right) $ as above, and 
\begin{equation*}
f_{n}\left( x,y\right) =\left( nv_{n}\right) ^{-1/2}\sqrt{q\left( \frac{x}{n}%
\right) q\left( \frac{y}{n}\right) };\;\;g_{n}\left( x\right) =\frac{%
n^{-1}\left( e^{ix}-1\right) }{e^{ix/n}-1}.
\end{equation*}%
The function $g$ is bounded by $2h\left( x\right) $. Elementary calculations
with $n\geq 3$ and $x,y\in I_{n}$ lead to 
\begin{equation}
\left\vert g_{n}\left( x+y\right) -g\left( x+y\right) \right\vert \leq \frac{%
\left\vert g\left( x+y\right) \right\vert }{1+n/\left\vert x+y\right\vert }%
\leq \frac{2h\left( x+y\right) }{1+n/\left\vert x+y\right\vert }.
\label{gn-g}
\end{equation}%
From Theorem \ref{T2chaos},%
\begin{eqnarray*}
&&d_{TV}\left( F_{n},F_{\infty }\right) \\
&\leq &\iint_{\mathbf{R}^{2}\backslash \left( I_{n}\right) ^{2}}\left\vert
f\left( x,y\right) \right\vert ^{2}h^{2}\left( x+y\right) dxdy \\
&&+\iint_{\left( I_{n}\right) ^{2}}\left( \left\vert f_{n}\left( x,y\right)
-f\left( x,y\right) \right\vert ^{2}\left\vert g_{n}\left( x+y\right)
\right\vert ^{2}+\left\vert f\left( x,y\right) \right\vert ^{2}\left\vert
g_{n}\left( x+y\right) -g\left( x+y\right) \right\vert ^{2}\right) dxdy \\
&\leq &\iint_{\mathbf{R}^{2}\backslash \left( I_{n}\right) ^{2}}\left\vert
f\left( x,y\right) \right\vert ^{2}h^{2}\left( x+y\right)
dxdy+16\iint_{\left( I_{n}\right) ^{2}}\left\vert f_{n}\left( x,y\right)
-f\left( x,y\right) \right\vert ^{2}h^{2}\left( x+y\right) dxdy \\
&&+\iint_{\left( I_{n}\right) ^{2}}\left\vert g_{n}\left( x+y\right)
-g\left( x+y\right) \right\vert ^{2}\left\vert f\left( x,y\right)
\right\vert ^{2}dxdy.
\end{eqnarray*}%
All three terms above converge to $0$, the first two by assumption, the last
by dominated convergence. However, to derive quantitative estimates, it is
best to exploit relation (\ref{gn-g}) more specifically. Therefore, for $%
\alpha \in \left( 0,1\right) $ fixed, we write%
\begin{eqnarray*}
&&\iint_{\left( I_{n}\right) ^{2}}\left\vert g_{n}\left( x,y\right) -g\left(
x,y\right) \right\vert ^{2}\left\vert f\left( x,y\right) \right\vert ^{2}dxdy
\\
&=&\left( \iint_{\left( I_{n}\right) ^{2}\backslash \left( I_{n^{\alpha
}}\right) ^{2}}+\iint_{\left( I_{n^{\alpha }}\right) ^{2}}\right) \left\vert
g_{n}\left( x,y\right) -g\left( x,y\right) \right\vert ^{2}\left\vert
f\left( x,y\right) \right\vert ^{2}dxdy \\
&\leq &4\iint_{\mathbf{R}^{2}\backslash \left( I_{n^{\alpha }}\right)
^{2}}\left\vert f\left( x,y\right) \right\vert ^{2}h^{2}\left( x+y\right)
dxdy+\frac{4}{n^{2-2\alpha }}\iint_{\mathbf{R}^{2}}\left\vert f\left(
x,y\right) \right\vert ^{2}h^{2}\left( x+y\right) dxdy.
\end{eqnarray*}%
The corollary easily follows.\hfill $\blacksquare $

\section{Example: log-modulation and second-chaos limits\label{LOG2}}

We have in mind the same class of examples as in Section \ref{LOG}. If $H>%
\frac{3}{4}$, one gets, by relation (\ref{equiv3}) in Theorem \ref{theo-est}%
, that $\kappa _{3}(F_{n})$ converges to a non-zero constant. Theorem \ref%
{T3M} then proves that $F_{n}$ cannot converge to a normal law. We now study
the possible convergence of $F_{n}$ to non-normal laws.

In order to streamline the presentation, since our strategy is to use
Corollary \ref{Cor2chaos}, it turns out to be more convenient to make
assumptions on $q$ and derive corresponding estimates on $\rho $ and other
quantities of interest. The reader may see below, by comparing our
definition of $q$ in (\ref{qyoo}) and the estimate on $\rho $ in Proposition %
\ref{prop-rho}, that assumptions on $q$ and $\rho $ of log-modulated power
type are asymptotically equivalent.

Let $\beta \geq 0$ and $H\in (3/4,1)$. Consider the positive function $q$ on
the unit circle $S^{1}=[-\pi ,\pi ]$ defined, except at $x=0$, by%
\begin{equation}
q\left( x\right) :=C_{H,\beta }\left\vert x\right\vert ^{1-2H}\log ^{2\beta
}\left( \frac{e\pi }{\left\vert x\right\vert }\right) .  \label{qyoo}
\end{equation}%
The constant $C_{H,\beta }$ is chosen in such a way that $\int_{-\pi }^{\pi
}q\left( x\right) dx/\left( 2\pi \right) =1$, i.e. 
\begin{equation}
C_{H,\beta }:=2\pi /\int_{-\pi }^{\pi }\left\vert x\right\vert ^{1-2H}\log
^{2\beta }\left( \frac{e\pi }{\left\vert x\right\vert }\right) dx,
\label{CHbeta}
\end{equation}%
in order to stay with the assumption that our stationary sequence has unit
variance, but other normalizing constants pose no additional difficulty. The
case $\beta =0$ corresponds to fGn. The case $\beta <0$ has slightly
different properties than the case $\beta \geq 0$, and requires further
computations; we omit it for the sake of conciseness. To simplify the
notation, we introduce 
\begin{equation*}
L\left( y\right) :=\log ^{2\beta }\left( \left\vert y\right\vert \right)
\end{equation*}%
and notice that $L\left( e\pi /\left\vert x\right\vert \right) \geq 1$ for
all $x\in \lbrack -\pi ,\pi ]$. Moreover, $q$ is in $L^{1}\left(
S^{1},dx\right) $ and is $C^{\infty }$ everywhere except at $0$, and
therefore $q$ coincides with the Fourier series of its Fourier inverse $\rho 
$; in other words, the stationary Gaussian process $X$ with covariance
function $\rho $ given by%
\begin{equation}
\rho \left( k\right) :=\frac{1}{2\pi }\int_{-\pi }^{\pi }q\left( x\right)
\cos \left( kx\right) dx  \label{rhodefq}
\end{equation}%
has spectral density $q$. The relation (\ref{rhodefq}), which is our
definition of $\rho $, serves as the Fourier inversion property required for
applying Theorem \ref{T2chaos} and its corollary. The other needed
assumption is $q\in L^{1}\left( S^{1}\right) $, which holds since $H<1$. We
recompute $\rho $ by changing variables and using the definition of $q$, to
get%
\begin{equation}
\rho \left( k\right) =\frac{C_{H,\beta }}{2\pi }k^{2H-2}\int_{-k\pi }^{k\pi
}\left\vert x\right\vert ^{1-2H}\cos \left( x\right) L\left( k\frac{e\pi }{%
\left\vert x\right\vert }\right) dx.  \label{rhonew}
\end{equation}%
Here normalizations are such that $\rho \left( 0\right) =1$. The asymptotics
of $\rho $ are the following.

\begin{prop}
\label{prop-rho}With $\rho $ in (\ref{rhonew}), with $C_{H,\beta }$ in (\ref%
{CHbeta}) and with%
\begin{equation*}
K_{H}:=\frac{1}{\pi }\int_{0}^{\infty }\left\vert x\right\vert ^{1-2H}\cos
\left( x\right) dx=2\Gamma \left( 2-2H\right) \cos \left( \pi \left(
1-H\right) \right) ,
\end{equation*}%
we have for large $k$%
\begin{equation*}
\rho \left( k\right) =C_{H,\beta }K_{H}L\left( k\right) k^{2H-2}\left( 1+%
\mathcal{O}\left( \frac{1}{L\left( k\right) }\right) \right) .
\end{equation*}
\end{prop}

The proof of this proposition is given in Section \ref{Aprop-rho} in the
Appendix. The next proposition, which gives the behavior of $%
nv_{n}=Var\left( \sum_{k=0}^{n-1}X_{k}^{2}\right) $ using Proposition \ref%
{prop-rho}, is also proved in the Appendix, in Section \ref{APnvn}. Recall
from (\ref{equivv}) in Theorem \ref{theo-est} that $nv_{n}=n+2%
\sum_{k=1}^{n-1}\left( n-k\right) \rho ^{2}\left( k\right) $.

\begin{prop}
\label{Pnvn}With $\rho $ in (\ref{rhonew}) and $nv_{n}=n+2\sum_{k=1}^{n-1}%
\left( n-k\right) \rho ^{2}\left( k\right) $, for $n$ large%
\begin{equation*}
nv_{n}=\left( C_{H,\beta }\right) ^{2}K_{H}^{\prime }~n^{4H-2}~L^{2}\left(
n\right) \left( 1+\mathcal{O}\left( \frac{1}{\log n}\right) \right)
\end{equation*}%
where $C_{H,\beta }$ is given in (\ref{CHbeta}) and 
\begin{equation*}
K_{H}^{\prime }:=\frac{\left( 2\Gamma \left( 2-2H\right) \cos \left( \pi
\left( 1-H\right) \right) \right) ^{2}}{\left( 4H-2\right) \left(
4H-3\right) }.
\end{equation*}
\end{prop}

With Proposition \ref{Pnvn} in hand, we can use Corollary \ref{Cor2chaos} to
establish a speed of convergence result for the Dobrushin-Major theorem.

\begin{theo}
\label{TNNS}Let $(F_{n})_{n\geq 0}$ be the sequence of normalized quadratic
variations of a centered stationary Gaussian sequence $\left( X_{n}\right) $
with unit variance, i.e. $F_{n}:=V_{n}/\sqrt{v_{n}}$, where $%
V_{n}:=n^{-1/2}\sum_{k=0}^{n-1}(X_{k}^{2}-1)$, and $%
v_{n}:=E[V_{n}^{2}]=n^{-1}\sum_{k=0}^{n-1}\sum_{\ell =0}^{n-1}\rho \left(
k-\ell \right) ^{2}$, with the covariance $\rho $ of $X$ given by (\ref%
{rhonew}). Assume $H\in (3/4,1)$ and $\beta \geq 0$. Then $F_{n}$ converges
in law to the law of a second-chaos variable $F_{\infty }$, and we have the
representation%
\begin{equation*}
F_{\infty }=\iint_{\mathbf{R}^{2}}\frac{\left\vert xy\right\vert ^{H-1/2}}{%
\sqrt{K_{H}^{\prime }}}\frac{e^{i(x+y)}-1}{i\left( x+y\right) }W\left(
dx\right) W\left( dy\right)
\end{equation*}%
where $W$ is a complex-valued white noise on $\mathbf{R}$ with scale
determined by $\left\vert W\left( dx\right) \right\vert ^{2}=dx/\left( 2\pi
\right) $ and $K_{H}^{\prime }$ is given in Proposition \ref{Pnvn}.
Moreover, the speed of convergence in total variation is bounded above as%
\begin{equation*}
d_{TV}\left( F_{n},F_{\infty }\right) \leq \frac{c}{\log ^{1/2}n}
\end{equation*}%
where the positive constant $c$ depends only on $H$ and $\beta $.
\end{theo}

Before concluding our article with the proof of this theorem, we rephrase
the result to match the usual normalizations found in the literature when $%
H>3/4$.

\begin{rem}
According to Theorem \ref{TNNS}, with $\left( X_{n}\right) $ the correlated
sequence of standard normal random variables defined therein, the sequence%
\begin{equation*}
\frac{\sum_{k=0}^{n-1}(X_{k}^{2}-1)}{n^{2H-1}}
\end{equation*}%
converges in distribution to the law of the second-chaos variable
represented by%
\begin{equation*}
C_{H,\beta }\iint_{\mathbf{R}^{2}}\left\vert xy\right\vert ^{H-1/2}\frac{%
e^{i(x+y)}-1}{i\left( x+y\right) }W\left( dx\right) W\left( dy\right)
\end{equation*}%
where $W$ is a complex-valued white noise on $\mathbf{R}$ with scale
determined by $\left\vert W\left( dx\right) \right\vert ^{2}=dx/\left( 2\pi
\right) $ and $C_{H,\beta }$ is given in (\ref{CHbeta}). The speed of
convergence in total variation is of order $\log ^{-1/2}n$ for any $\beta >0$
and any $H\in (3/4,1)$. When $\beta =0$, this speed can be improved to $%
n^{3/4-H}$ as can be seen in \cite[Theorem 7.4.5]{npBook}.
\end{rem}

\noindent \emph{Proof of Theorem \ref{TNNS}.} Since, as mentioned above, $q$
and $\rho $ satisfy the assumptions of Theorem \ref{T2chaos}, according to
Corollary \ref{Cor2chaos}, and using its notation, it is sufficient to prove
that, with $f\left( x,y\right) =\left\vert xy\right\vert ^{H-1/2}/\sqrt{%
K_{H}^{\prime }},$ the three terms in lines (\ref{cor2chaosline1}), (\ref%
{cor2chaosline2}), and (\ref{cor2chaosline3}) all converge to 0 at least as
fast as $\log ^{-2}n$. In all calculations below, it will be convenient to
change variables and replace $y$ by $-y$; this only entails replacing $%
h\left( x+y\right) $ by $h\left( x-y\right) $; all other expressions are
invariant by this change of variables since $I_{n}$ and $q$ are even and $f$
is even in $y$; see Remark \ref{rem-even}.\vspace*{0.1in}

\noindent \emph{Step 1}. We use Proposition \ref{Pnvn} to handle the term in
line (\ref{cor2chaosline2}). With $\gamma =2-2H$, this term is bounded above
by a constant times%
\begin{eqnarray*}
\mathcal{I}\left( 2\right) &:&=\iint_{\left( I_{n}\right) ^{2}}\left\vert 
\sqrt{\frac{q\left( \frac{x}{n}\right) q\left( \frac{y}{n}\right) }{nv_{n}}}%
-f\left( x,y\right) \right\vert ^{2}h^{2}\left( x-y\right) dxdy \\
&=&\iint_{\left( I_{n}\right) ^{2}}\frac{\left\vert xy\right\vert ^{\gamma
-1}}{K_{H}^{\prime }}\left\vert \frac{\sqrt{L\left( \frac{e\pi n}{x}\right)
L\left( \frac{e\pi n}{y}\right) }}{\left( 1+\mathcal{O}\left( \frac{1}{\log n%
}\right) \right) L\left( n\right) }-1\right\vert ^{2}\min \left(
1,\left\vert x-y\right\vert ^{-2}\right) dxdy.
\end{eqnarray*}%
We can write 
\begin{equation*}
\frac{L\left( \frac{e\pi n}{x}\right) }{L\left( n\right) }=\left( 1-\frac{%
\log \left\vert x\right\vert }{\log \left( e\pi n\right) }\right) ^{2\beta }.
\end{equation*}%
In order to use a uniform bound on the first order Taylor expansion for this
expression, we will fix $\delta \in \left( 0,1\right) $ and consider $u\in
\lbrack -\delta ,1]$; then there is a positive finite constant $c=c\left(
\beta ,\delta \right) $ such that for large $n$, $\left\vert \left(
1-u\right) ^{\beta }-1+\beta u\right\vert \leq cu^{2}$. Hence for $%
\left\vert x\right\vert \in \lbrack n^{-\delta },\pi n]$ and $n$ large
enough,%
\begin{equation*}
\frac{\sqrt{L\left( \frac{e\pi n}{x}\right) L\left( \frac{e\pi n}{y}\right) }%
}{\left( 1+\mathcal{O}\left( \frac{1}{\log n}\right) \right) L\left(
n\right) }=1+\mathcal{O}\left( \frac{1}{\log n}\right) -\beta \frac{\log
\left\vert x\right\vert +\log \left\vert y\right\vert }{\log n}+c\frac{\log
^{2}\left\vert x\right\vert +\log ^{2}\left\vert y\right\vert }{\log ^{2}n}.
\end{equation*}%
\vspace*{0.1in}

\noindent \emph{Step 1.1}. For $n$ large enough, by symmetry, the portion of 
$\mathcal{I}\left( 2\right) $ for $\left\vert x\right\vert \geq n^{-\delta }$
and $\left\vert y\right\vert \geq n^{-\delta }$ is bounded above by a
constant times%
\begin{equation*}
\int_{0}^{\pi n}\int_{0}^{x}\left( xy\right) ^{\gamma -1}\left( \frac{1}{%
\log n}+\frac{\left\vert \log x\right\vert +\left\vert \log y\right\vert }{%
\log n}+\frac{\log ^{2}x+\log ^{2}y}{\log ^{2}n}\right) ^{2}\min \left(
1,\left\vert x-y\right\vert ^{-2}\right) dydx.
\end{equation*}%
Since $\left( xy\right) ^{\gamma -1}\left( 1+\left\vert \log x\right\vert
+\log ^{2}x\right) \left( 1+\left\vert \log y\right\vert +\log ^{2}y\right) 
\mathbf{1}_{0\leq y\leq x}$ is integrable for $x$ near $0$, the contribution
of the above integral for $x\in \lbrack 0,2]$ is $\mathcal{O}\left( \log
^{-2}n\right) $. Therefore, by separating this set and the diagonal set $%
x-1<y<x$ from the rest, the above integral is bounded above by a constant
times%
\begin{eqnarray*}
&&\frac{1}{\log ^{2}n}+\int_{2}^{\pi n}\int_{x-1}^{x}y^{\gamma -1}\frac{\log
^{2}x+\log ^{2}y+\log ^{4}x+\log ^{4}y}{\log ^{2}n}dydx \\
&&+\int_{2}^{\pi n}\int_{0}^{x-1}y^{\gamma -1}\left( x-y\right) ^{-2}\frac{%
\log ^{2}x+\log ^{2}y+\log ^{4}x+\log ^{4}y}{\log ^{2}n}dydx.
\end{eqnarray*}%
This is itself bounded above by%
\begin{eqnarray*}
&&\frac{1}{\log ^{2}n}+\frac{4}{\log ^{2}n}\int_{2}^{\pi n}\left( x-1\right)
^{\gamma -1}\log ^{4}x~dx \\
&&+\frac{1}{\log ^{2}n}\int_{2}^{\pi n}\int_{0}^{1}y^{\gamma -1}\left(
x-y\right) ^{-2}\left( \log ^{2}\left( x/y\right) +\log ^{4}\left(
x/y\right) \right) dydx \\
&&+\frac{4}{\log ^{2}n}\int_{2}^{\pi n}\log ^{4}x\int_{1}^{x-1}y^{\gamma
-1}\left( x-y\right) ^{-2}dydx \\
&\leq &\frac{1}{\log ^{2}n}+\frac{4}{\log ^{2}n}\left( \int_{1}^{\infty
}x^{\gamma -1}\log ^{4}\left( x+1\right) dx\right) \\
&&+\frac{1}{\log ^{2}n}\int_{2}^{\pi n}x^{\gamma -2}\int_{0}^{1/x}z^{\gamma
-1}\left( 1-z\right) ^{-2}\left( \log ^{2}z+\log ^{4}z\right) dzdx \\
&&+\frac{4}{\log ^{2}n}\int_{2}^{\pi n}x^{\gamma -2}\log
^{4}x\int_{1/x}^{1-1/x}z^{\gamma -1}\left( 1-z\right) ^{-2}dydx \\
&\leq &\frac{c}{\log ^{2}n}\left( \int_{2}^{\pi n}x^{\gamma -2}x^{-\gamma
+\varepsilon }dx+\int_{2}^{\pi n}x^{\gamma -2}\log ^{4}x dx\left[
4\int_{1/x}^{1/2}z^{\gamma -1}dz+4\int_{1/x}^{1/2}z^{-2}dz\right] \right) \\
&\leq &\frac{c}{\log ^{2}n}\left( 1+\int_{2}^{\pi n}x^{\gamma -2}\left(
x^{-\gamma +\varepsilon }+x^{-\gamma }+x^{-1}\right) \log ^{4}x dx\right) ,
\end{eqnarray*}%
for any constant $\varepsilon >0$, where $c$ depends only on $\gamma $ and $%
\varepsilon $. By choosing $\varepsilon <\gamma $, we get that this portion
of $\mathcal{I}\left( 2\right) $ goes to $0$ as fast as $\log ^{-2}n$.%
\vspace*{0.1in}

\noindent \emph{Step 1.2}. For $n$ large enough, by symmetry, the portion of 
$\mathcal{I}\left( 2\right) $ for $\left\vert x\right\vert \leq n^{-\delta }$
or $\left\vert y\right\vert \leq n^{-\delta }$ is bounded above by a
constant times%
\begin{eqnarray*}
&&\int_{0}^{\pi n}y^{\gamma -1}dy\int_{0}^{n^{-\delta }}x^{\gamma -1}\left(
1+\frac{L\left( \frac{e\pi n}{x}\right) }{L\left( n\right) }+\frac{L\left( 
\frac{e\pi n}{y}\right) }{L\left( n\right) }\right) \min \left(
1,\left\vert x-y\right\vert ^{-2}\right) dx \\
&\leq &\int_{0}^{n^{-\delta }}y^{\gamma -1}dy\int_{0}^{n^{-\delta
}}x^{\gamma -1}\left( 2+L\left( \frac{1}{x}\right) +L\left( \frac{1}{y}%
\right) \right) dx \\
&&+\int_{n^{-\delta }}^{1}y^{\gamma -1}dy\int_{0}^{n^{-\delta }}x^{\gamma
-1}\left( 3+L\left( \frac{1}{x}\right) \right) dx \\
&&+\int_{n^{-\delta }}^{\pi n}y^{\gamma -1}dy\int_{0}^{n^{-\delta
}}x^{\gamma -1}\left( 3+L\left( \frac{1}{x}\right) \right) dx\frac{4}{y^{2}}.
\end{eqnarray*}%
Therefore, for any $\varepsilon >0$, and $n$ large enough, the above
expression is bounded above by a constant times%
\begin{eqnarray*}
&&\int_{0}^{n^{-\delta }}y^{\gamma -1}dy~n^{-\gamma \delta }\left( 2+L\left( 
\frac{1}{y}\right) +n^{\varepsilon }\right) +\int_{n^{-\delta
}}^{1}y^{\gamma -1}dy~n^{-\gamma \delta }\left( 3+n^{\varepsilon }\right)
+\int_{1}^{\pi n}y^{\gamma -3}dy~n^{-\gamma \delta }\left( 3+n^{\varepsilon
}\right) \\
&\leq &\gamma ^{-1}n^{-\gamma \delta }\left( 2+n^{\varepsilon
}+n^{2\varepsilon }\right) +4\gamma ^{-1}n^{-\gamma \delta +\varepsilon
}+4\left( 2-\gamma \right) ^{-1}n^{-\gamma \delta +\varepsilon } \\
&\leq &cn^{-\gamma \delta +2\varepsilon }
\end{eqnarray*}%
for some constant $c$ depending only on $\gamma $ and $\varepsilon $. For $%
\varepsilon <\gamma \delta /2$, we get that this portion of $\mathcal{I}%
\left( 2\right) $ goes to $0$ faster than $\log ^{-2}n$. Hence $\mathcal{I}%
\left( 2\right) =\mathcal{O}\left( \log ^{-2}n\right) $. \vspace*{0.1in} 
\newline

\noindent \emph{Step 2}. We now consider the term in line (\ref%
{cor2chaosline1}). By symmetry, this is bounded above by a constant times%
\begin{equation*}
\mathcal{I}\left( 1\right) :=\int_{0}^{\infty }x^{\gamma -1}\int_{n^{\alpha
}}^{\infty }y^{\gamma -1}\min \left( 1,\left\vert x-y\right\vert
^{-2}\right) dydx.
\end{equation*}%
That is in turn bounded above by%
\begin{equation*}
\int_{0}^{n^{\alpha }-1}x^{\gamma -1}\left( \int_{n^{\alpha }}^{\infty
}y^{\gamma -1}\left( y-x\right) ^{-2}dy\right) dx+2\int_{n^{\alpha
}-1}^{\infty }x^{\gamma -1}\int_{n^{\alpha }-1}^{x}y^{\gamma -1}\min \left(
1,\left\vert x-y\right\vert ^{-2}\right) dydx
\end{equation*}%
which is itself bounded above by%
\begin{eqnarray*}
&&n^{\alpha \left( \gamma -1\right) }\int_{0}^{n^{\alpha }-1}x^{\gamma
-1}\left( n^{\alpha }-x\right) ^{-1}dx+2\int_{n^{\alpha }-1}^{\infty
}x^{\gamma -1}\int_{n^{\alpha }-1}^{x-1}y^{\gamma -1}\left( x-y\right)
^{-2}dydx \\
&&+2\int_{n^{\alpha }-1}^{\infty }x^{\gamma -1}\int_{x-1}^{x}y^{\gamma
-1}dydx \\
&\leq &n^{\alpha \left( \gamma -1\right) }\left( \frac{2}{n^{\alpha }}%
\int_{0}^{n^{\alpha }/2}x^{\gamma -1}dx+\left( \frac{2}{n^{\alpha }}\right)
^{\gamma -1}\int_{1}^{n^{\alpha }/2}x^{-1}dx\right) \\
&&+2\int_{n^{\alpha }-1}^{\infty }y^{\gamma -1}\int_{x+1}^{\infty }\left(
x-y\right) ^{-2}dxdy+2\int_{n^{\alpha }-1}^{\infty }\left( x-1\right)
^{2\gamma -2}dx \\
&=&n^{2\alpha \left( \gamma -1\right) }\left( \frac{2^{1-\gamma }}{\gamma }%
+2^{\gamma -1}\log \left( \frac{n^{\alpha }}{2}\right) \right) +\frac{2}{%
1-\gamma }\left( n^{\alpha }-1\right) ^{\gamma -1}+\frac{2}{1-2\gamma }%
\left( n^{\alpha }-2\right) ^{2\gamma -1}.
\end{eqnarray*}%
Since $\gamma =2-2H$ and $H\in (3/4,1)$, we have $0>2\gamma -1>\gamma
-1>2\gamma -2,$ and thus the last expression above is a $\mathcal{O}\left(
n^{-\alpha \left( 1-2\gamma \right) }\right) $ which goes to $0$ faster than 
$\log ^{-2}n$. I.e. $\mathcal{I}\left( 1\right) =\mathcal{O}\left(
n^{-\alpha \left( 1-2\gamma \right) }\right) \ll \log ^{-2}n$.\vspace*{0.1in}

\noindent \emph{Step 3}. Finally we consider the term in line (\ref%
{cor2chaosline3}). By symmetry, this is bounded above by a constant times%
\begin{equation*}
\mathcal{I}\left( 3\right) :=\frac{1}{n^{2-2\alpha }}\int_{0}^{\infty
}x^{\gamma -1}\int_{0}^{x}y^{\gamma -1}\min \left( 1,\left( x-y\right)
^{-2}\right) dydx.
\end{equation*}%
To show that this $\mathcal{I}\left( 3\right) =o\left( \log ^{-2}n\right) $
it is sufficient to prove that the integral above is finite. Since $\gamma
-1>-1$, The portion of that integral corresponding to $x\in \lbrack 0,2]$ is
finite. We thus only need to study the portion corresponding to $x>2$ :%
\begin{eqnarray*}
&&\int_{2}^{\infty }x^{\gamma -1}\int_{0}^{x}y^{\gamma -1}\min \left(
1,\left( x-y\right) ^{-2}\right) dydx \\
&\leq &\int_{2}^{\infty }x^{\gamma -1}\left( \left( \frac{x}{2}\right)
^{-2}\int_{0}^{x/2}y^{\gamma -1}dy+\left( \frac{x}{2}\right) ^{\gamma
-1}\int_{x/2}^{x-1}\left( x-y\right) ^{-2}dy+\int_{x-1}^{x}y^{\gamma
-1}dy\right) dx \\
&\leq &\int_{2}^{\infty }x^{\gamma -1}\left( \left( \frac{x}{2}\right)
^{-2}\gamma ^{-1}\left( \frac{x}{2}\right) ^{\gamma }+\left( \frac{x}{2}%
\right) ^{\gamma -1}\left( \int_{1}^{\infty }y^{-2}dy\right) +\left(
x-1\right) ^{\gamma -1}\right) dx \\
&\leq &\gamma ^{-1}2^{2-\gamma }\int_{2}^{\infty }x^{2\gamma
-3}dx+2^{2-\gamma }\int_{1}^{\infty }x^{2\gamma -2}dx.
\end{eqnarray*}%
This is finite, given that $H>3/4\implies 2\gamma -2<-1$. The proof of the
theorem is complete. $\hspace*{\fill}\blacksquare $\bigskip

\section{Appendix: technical elements used in Sections \protect\ref{SNN} and 
\protect\ref{LOG2}.}

\subsection{Representations of stationary Gaussian processes using the
complex-valued white noise measure\label{Complex}}

For the reader's convenience, and for the sake of being self-contained to
some extent, we briefly recall the construction and properties of the
complex-valued white noise $W$. As indicated in Section \ref{SNN}, $W$ is
the independently scattered $\mathbf{C}$-valued centered Gaussian measure on 
$S^{1}=[-\pi ,\pi ]$ such that for $x\in \lbrack 0,\pi ]$, $W\left(
dx\right) =B_{1}\left( dx\right) +iB_{2}\left( dx\right) $, where $B_{1}$
and $B_{2}$ are two real-valued independent white noise measures on $[0,\pi
] $ such that $Var[B_{i}([0,\pi ])]=1/2$, and for every $x\in \lbrack 0,\pi
] $, $W\left( -dx\right) =\overline{W\left( dx\right) }=B_{1}\left(
dx\right) -iB_{2}\left( dx\right) $.

This definition of $W$, where one notes that unlike in the real case, the
restrictions of $W$ to $[-\pi ,0]$ and to $[0,\pi ]$ are \emph{not}
independent, implies the following properties, using the usual shorthand
differential notation for It\^{o}'s rule: for all $x\in \lbrack 0,\pi ]$,

\begin{itemize}
\item first It\^{o} rule: $W\left( dx\right) W\left( -dx\right) =W\left(
dx\right) \overline{W\left( dx\right) }=\left\vert W\left( dx\right)
\right\vert ^{2}=dx/\left( 2\pi \right) ,$

\item second It\^{o} rule: $W\left( dx\right) W\left( dx\right) =W\left(
dx\right) ^{2}=0,$

\item $W\left( dx\right) $ and $W\left( dy\right) $ are independent for $%
x\neq y$ and $xy\geq 0$,
\end{itemize}

For $X$ defined as in (\ref{SR}) using a non-negative function $q\in
L^{1}\left( S^{1}\right) $, i.e. for $k\in \mathbf{Z}$,%
\begin{equation*}
X(k)=\int_{S^{1}}e^{ikx}\sqrt{q(x)}W(dx),
\end{equation*}%
we can rewrite this expression by expanding $W$ according to its definition
above: for all $k\in \mathbf{Z}$,%
\begin{eqnarray}
X\left( k\right)  &=&2\int_{0}^{\pi }\sqrt{q\left( x\right) }\cos \left(
kx\right) B_{1}\left( dx\right) -2\int_{0}^{\pi }\sqrt{q\left( x\right) }%
\sin \left( kx\right) B_{2}\left( dx\right)   \label{SR2} \\
&=&\int_{-\pi }^{\pi }\sqrt{q\left( x\right) }\cos \left( kx\right)
W_{1}\left( dx\right) +\int_{-\pi }^{\pi }\sqrt{q\left( x\right) }\sin
\left( kx\right) W_{2}\left( dx\right)   \label{SRusual}
\end{eqnarray}%
where the second equality is in law, with $W_{1}$ and $W_{2}$ two
independent real-valued white noise measures on $[-\pi ,\pi ]$ with scale
determined by $Var\left[ W_{i}\left( [-\pi ,\pi ]\right) \right] =1$.

Representation (\ref{SRusual}) is more commonly found in the literature than
Representation (\ref{SR2}). From either representation, it is evident that $%
X $ is real valued. From (\ref{SR}), one can check that $X$ has the
announced covariance, as follows. Using the fact that $X$ is real-valued,
and the product rule for Wiener integrals, we have $X\left( k\right) X\left(
l\right) =X\left( k\right) \overline{X\left( l\right) }$, which is the sum
of a mean-zero second-chaos variable and of a constant formally expressed as 
$\int_{S^{1}}e^{ikx}\sqrt{q(x)}W(dx)e^{-ilx}\sqrt{q(x)}\ \overline{W\left(
dx\right) }$. Thus, by It\^{o}'s rule,%
\begin{equation*}
\mathbf{E}\left[ X\left( k\right) X\left( l\right) \right] =\int_{-\pi
}^{\pi }e^{i\left( k-l\right) x}q(x)\frac{dx}{2\pi }=\rho \left( k-l\right) ,
\end{equation*}%
where the last equality assumes that Fourier inversion holds for $\rho $. We
can extend this type of calculation in general to express the isometry
property in the second chaos of the complex-valued $W$. For a function $f\in
L^{2}\left( [0,\pi ]^{2},\mathbf{C}\right) $ we can write%
\begin{eqnarray*}
I_{2}\left( f\right) &=&\iint_{[0,\pi ]^{2}}f\left( x,y\right) W\left(
dx\right) W\left( dy\right) +\iint_{[0,\pi ]\times \lbrack -\pi ,0]}f\left(
x,y\right) W\left( dx\right) W\left( dy\right) \\
&&+\iint_{[-\pi ,0]^{2}}f\left( x,y\right) W\left( dx\right) W\left(
dy\right) +\iint_{[-\pi ,0]\times \lbrack 0,\pi ]}f\left( x,y\right) W\left(
dx\right) W\left( dy\right) \\
&=&\iint_{[0,\pi ]^{2}}f\left( x,y\right) W\left( dx\right) W\left(
dy\right) +\iint_{[0,\pi ]^{2}}f\left( -x,-y\right) \overline{W\left(
dx\right) W\left( dy\right) } \\
&&+\iint_{[0,\pi ]^{2}}f\left( x,-y\right) W\left( dx\right) \overline{%
W\left( dy\right) }+\iint_{[0,\pi ]^{2}}f\left( -x,y\right) \overline{%
W\left( dx\right) }W\left( dy\right) .
\end{eqnarray*}%
When squaring $I_{2}\left( f\right) $ and taking its expectation, by the
second It\^{o} rule above, the only terms that remain are those for which a
product of the form $W\left( dx\right) W\left( dx\right) $ or $W\left(
dy\right) W\left( dy\right) $ does not appear. Thus we get only two terms
left:%
\begin{equation*}
\mathbf{E}\left[ I_{2}\left( f\right) ^{2}\right] =2\iint_{[0,\pi
]^{2}}f\left( x,y\right) f\left( -x,-y\right) \frac{dxdy}{\left( 2\pi
\right) ^{2}}+2\iint_{[0,\pi ]^{2}}f\left( x,-y\right) f\left( -x,y\right) 
\frac{dxdy}{\left( 2\pi \right) ^{2}}.
\end{equation*}%
When $f$ satisfies the Hermitian evenness property $f\left( -x,-y\right) =%
\overline{f\left( x,y\right) }$, this formula easily yields the isometry
property (\ref{isometry}) announced earlier:%
\begin{equation*}
\mathbf{E}\left[ I_{2}\left( f\right) ^{2}\right] =\iint_{[-\pi ,\pi
]^{2}}\left\vert f\left( x,y\right) \right\vert ^{2}\frac{dxdy}{\left( 2\pi
\right) ^{2}}=:\left\Vert f\right\Vert _{L^{2}\left( \left( S^{1}\right)
^{2}\right) }^{2}.
\end{equation*}

\subsection{Proof of Proposition \protect\ref{prop-rho}\label{Aprop-rho}}

The constant 
\begin{equation*}
\int_{\mathbf{R}}\left\vert x\right\vert ^{1-2H}\cos \left( x\right) L\left( 
\frac{e\pi }{\left\vert x\right\vert }\right) dx
\end{equation*}%
is finite because $x\mapsto \left\vert x\right\vert ^{1-2H}L\left( e\pi
/\left\vert x\right\vert \right) $ decreases to $0$ as $x\rightarrow \infty $%
. Thus by (\ref{rhonew}), we get

\begin{equation*}
\frac{\rho \left( k\right) }{C_{H,\beta }}=\frac{k^{2H-2}}{\pi }%
\int_{0}^{\infty }\left\vert x\right\vert ^{1-2H}\cos \left( x\right)
L\left( k\frac{e\pi }{\left\vert x\right\vert }\right) dx-\frac{k^{2H-2}}{%
\pi }\int_{k\pi }^{\infty }\left\vert x\right\vert ^{1-2H}\cos \left(
x\right) dx
\end{equation*}%
After integrating by parts, the second term on the right-hand side above can
be written as 
\begin{equation*}
\frac{k^{2H-2}}{\pi }\lim_{N\rightarrow \infty }\int_{k\pi }^{N}x^{1-2H}\cos
\left( x\right) dx=\frac{k^{2H-2}}{\pi }\lim_{N\rightarrow \infty }\left(
N^{1-2H}\sin N+\left( 2H-1\right) \int_{k\pi }^{N}x^{-2H}\sin \left(
x\right) dx\right) =\mathcal{O}\left( k^{-1}\right) .
\end{equation*}%
Therefore to prove the proposition, it is sufficient to show that%
\begin{equation*}
\frac{1}{\pi }\int_{0}^{\infty }x^{1-2H}\cos \left( x\right) L\left( k\frac{%
e\pi }{x}\right) dx=K_{H}L\left( k\right) \left( 1+\mathcal{O}\left( \frac{1%
}{L\left( k\right) }\right) \right) .
\end{equation*}%
There is a positive constant $c\left( \beta \right) $ such that for $%
\left\vert y\right\vert <1$, we have $1-c\left( \beta \right) \left\vert
y\right\vert \leq \left( 1+y\right) ^{\beta }\leq 1+c\left( \beta \right)
\left\vert y\right\vert $. Thus for $0\leq y\leq 1$, we can write $\left(
1+y\right) ^{\beta }=1+\mathcal{O}\left( y\right) $ where $\left\vert 
\mathcal{O}\left( y\right) /y\right\vert $ is bounded by $c\left( \beta
\right) $. We also use $\left\vert \cos x\right\vert \leq 1$ when $x$ is
small. We compute, for any $\varepsilon >0$, for $k$ large, 
\begin{eqnarray*}
&&\int_{0}^{\infty }x^{1-2H}\cos \left( x\right) L\left( k\frac{e\pi }{x}%
\right) dx \\
&=&\int_{0}^{e\pi /k}x^{1-2H}\cos \left( x\right) L\left( k\frac{e\pi }{x}%
\right) dx +\int_{e\pi /k}^{\infty }\left\vert x\right\vert ^{1-2H}\cos \left(
x\right) \log ^{2\beta }k\left( 1+\frac{\log \left( \frac{e\pi }{\left\vert
x\right\vert }\right) }{\log k}\right) ^{2\beta }dx \\
&=&\mathcal{O}\left( \int_{0}^{e\pi /k}x^{1-2H}L\left( k\frac{e\pi }{x}%
\right) dx\right) +\int_{e\pi /k}^{\infty }\left\vert x\right\vert
^{1-2H}\cos \left( x\right) \log ^{2\beta }k\left( 1+\mathcal{O}\left( \frac{%
\log \left( \frac{e\pi }{\left\vert x\right\vert }\right) }{\log k}\right)
\right) dx \\
&=&\mathcal{O}\left( k^{2-2H}\int_{0}^{k^{-2}}y^{1-2H}\log ^{2\beta }\left(
y^{-1}\right) dy\right) +L\left( k\right) \int_{e\pi /k}^{\infty }\left\vert
x\right\vert ^{1-2H}\cos \left( x\right) dx \\
&&+\mathcal{O}\left( \frac{L\left( k\right) }{\log k}\right) \int_{e\pi
/k}^{\infty }\left\vert x\right\vert ^{1-2H}\cos \left( x\right) \log \left( 
\frac{e\pi }{\left\vert x\right\vert }\right) dx \\
&=&\mathcal{O}\left( k^{2-2H}k^{-2\left( 2-2H-\varepsilon \right) }\right)
+L\left( k\right) \pi K_{H}-L\left( k\right) \int_{0}^{e\pi /k}\left\vert
x\right\vert ^{1-2H}\cos \left( x\right) dx \\
&&+\mathcal{O}\left( \frac{L\left( k\right) }{\log k}\right)
\int_{0}^{\infty }\left\vert x\right\vert ^{1-2H}\cos \left( x\right) \log
\left( \frac{e\pi }{\left\vert x\right\vert }\right) dx \\
&=&\mathcal{O}\left( k^{-(2-2H)+\varepsilon }\right) +L\left( k\right) \pi
K_{H}+L\left( k\right) \mathcal{O}\left( k^{-(2-2H)}\right) +\mathcal{O}%
\left( \frac{L\left( k\right) }{\log k}\right) .
\end{eqnarray*}%
We used the facts that $y\mapsto \int_{0}^{y}x^{1-2H}dx=\mathcal{O}\left(
y^{2-2H}\right) $ near $0$, that for any $\varepsilon >0$, $y\mapsto
\int_{0}^{y}x^{1-2H}\log ^{\beta }\left( x^{-1}\right) dx=\mathcal{O}\left(
y^{2-2H-\varepsilon }\right) $ near $0$, and that $\int_{0}^{\infty
}\left\vert x\right\vert ^{1-2H}\cos \left( x\right) \log \left( \frac{e\pi 
}{\left\vert x\right\vert }\right) dx$ is a converging series. By taking $%
\varepsilon \in (0,2-2H)$, this proves the proposition. \hfill $\blacksquare 
$

\subsection{Proof of Proposition \protect\ref{Pnvn}\label{APnvn}}

Let $\alpha \in \left( 0,1\right) $. We compute $nv_{n}$ by splitting its
series up at the value $k=n^{\alpha }$. We use the notation $\gamma :=2-2H$
for compactness; note that $1>H>3/4$ implies $2\gamma \in \left( 0,1\right) $%
.

To lighten the notation slightly, we write $K_{H,\beta }:=K_{H}C_{H,\beta }$%
. From Proposition \ref{prop-rho}, we compute 
\begin{eqnarray}
nv_{n} &=&n+2K_{H,\beta }^{2}\left[ \sum_{k=1}^{[n^{\alpha
}]}+\sum_{k=[n^{\alpha }]+1}^{n}\right] \left( n-k\right) k^{-2\gamma
}L^{2}\left( k\right) \left( 1+\mathcal{O}\left( 1/L\left( k\right) \right)
\right) ^{2}  \notag \\
&=&n+2K_{H,\beta }^{2}\left( 1+\mathcal{O}\left( 1\right) \right)
n^{1+\alpha }\sum_{k=1}^{[n^{\alpha }]}k^{-2\gamma }L^{2}\left( k\right)
\left( 1+\mathcal{O}\left( 1/L\left( k\right) \right) \right) ^{2}
\label{linetwo} \\
&&+2K_{H,\beta }^{2}~n^{2-2\gamma }\frac{1}{n}\sum_{k=[n^{\alpha
}]+1}^{n}\left( 1-\frac{k}{n}\right) \left( \frac{k}{n}\right) ^{-2\gamma
}L^{2}\left( k\right) \left( 1+\mathcal{O}\left( 1/L\left( n^{\alpha
}\right) \right) \right) ^{2}.  \label{linethree}
\end{eqnarray}%
In line (\ref{linetwo}), the term $n$ is negligible in front of the
remainder of that line (we already knew this from Corollary \ref%
{Cor-for-nonnormal}), which is of order at least $n^{1+\alpha }\sum_{k=1}^{%
\left[ n^{\alpha }\right] }k^{-2\gamma }\asymp n^{1+\alpha +\alpha \left(
1-2\gamma \right) }=n^{1+\alpha \left( 1-2\gamma \right) }$ and no greater
than $n^{1+\alpha \left( 1-2\gamma \right) }L^{2}\left( n\right) $. Thus the
term in line (\ref{linetwo}) is $\mathcal{O}\left( n^{1+\alpha \left(
1-2\gamma \right) }L^{2}\left( n\right) \right) $. On the other hand, we set
up the term in line (\ref{linethree}) to draw a precise comparison with a
Riemann sum; thus modulo the factor $L^{2}\left( k\right) \left( 1+\mathcal{O%
}\left( 1/L\left( n^{\alpha }\right) \right) \right) ^{2}$ which is smaller
than any power, we have an expression which is asymptotically equivalent to
a constant times $n^{2-2\gamma }$. However we find $1+\alpha \left(
1-2\gamma \right) <2-2\gamma \iff 2\gamma <1$. Thus the terms in line (\ref%
{linethree}) dominate those in line (\ref{linetwo}) by a factor greater than
a small power. In other words, we have proved that for some $\varepsilon >0$,%
\begin{equation*}
nv_{n}=2K_{H,\beta }^{2}~\left( 1+\mathcal{O}\left( 1/L\left( n^{\alpha
}\right) \right) \right) ^{2}\left( 1+\mathcal{O}\left( n^{-\varepsilon
}\right) \right) ~n^{2-2\gamma }\frac{1}{n}\sum_{k=[n^{\alpha
}]+1}^{n}\left( 1-\frac{k}{n}\right) \left( \frac{k}{n}\right) ^{-2\gamma
}L^{2}\left( k\right) .
\end{equation*}%
Since $L\left( n^{\alpha }\right) =\mathcal{O}\left( L\left( n\right)
\right) $ trivially, and $\left( 1+\mathcal{O}\left( 1/L\left( n\right)
\right) \right) ^{p}=\left( 1+\mathcal{O}\left( 1/L\left( n\right) \right)
\right) $ for $p>0$, and for any integer $k\geq n^{\alpha }$ we can write%
\begin{equation*}
L\left( k\right) =\log ^{4\beta }\left( n\right) \left( 1+\frac{\log \left(
k/n\right) }{\log n}\right) ^{4\beta }=L\left( n\right) \left( 1+\mathcal{O}%
\left( \frac{1}{\log n}\right) \right) ,
\end{equation*}%
we get%
\begin{equation*}
nv_{n}=2K_{H,\beta }^{2}~\left( 1+\mathcal{O}\left( \frac{1}{L\left(
n\right) }\right) \right) ~n^{2-2\gamma }~L^{2}\left( n\right) ~\frac{1}{n}%
\sum_{k=[n^{\alpha }]+1}^{n}\left( 1-\frac{k}{n}\right) \left( \frac{k}{n}%
\right) ^{-2\gamma }.
\end{equation*}%
We must now compute the asymptotics of the series in the last line above.
Let $h_{\gamma }\left( x\right) :=\left( 1-x\right) x^{-2\gamma }$ defined
on $(0,1]$. We compute $h_{\gamma }^{\prime }\left( x\right) =x^{-2\gamma
}\left( 2\gamma -1-2\gamma /x\right) <0$ and find $\left\vert h_{\gamma
}^{\prime }\left( x\right) \right\vert \leq 2x^{-2\gamma -1}$. We thus have%
\begin{eqnarray*}
&&\left\vert \frac{1}{n}\sum_{k=[n^{\alpha }]+1}^{n}h_{\gamma }\left( \frac{k%
}{n}\right) -\int_{0}^{1}h_{\gamma }\left( x\right) dx\right\vert \leq \frac{%
1}{n}\int_{n^{\alpha -1}}^{1}2x^{-2\gamma -1}dx+\int_{0}^{n^{\alpha
-1}}x^{-2\gamma }dx \\
&=&\mathcal{O}\left( n^{2\gamma \left( 1-\alpha \right) -1}\right) +\mathcal{%
O}\left( n^{2\gamma \left( 1-\alpha \right) -1+\alpha }\right) =\mathcal{O}%
\left( n^{-\left( 1-\alpha \right) \left( 1-2\gamma \right) }\right) .
\end{eqnarray*}%
This proves that for any choice of $\alpha \in (0,1)$%
\begin{eqnarray*}
nv_{n} &=&2K_{H,\beta }^{2}~\left( 1+\mathcal{O}\left( \frac{1}{L\left(
n\right) }\right) \right) ~n^{2-2\gamma }~L^{2}\left( n\right)
~\int_{0}^{1}h_{\gamma }\left( x\right) dx~\left( 1+O\left( n^{-\left(
1-\alpha \right) \left( 1-2\gamma \right) }\right) \right) \\
&=&2K_{H,\beta }^{2}\left( \int_{0}^{1}h_{\gamma }\left( x\right) dx\right)
~n^{2-2\gamma }~L^{2}\left( n\right) ~\left( 1+\mathcal{O}\left( \frac{1}{%
L\left( n\right) }\right) \right)
\end{eqnarray*}%
which is the proposition's claim, given $\gamma =2-2H$ and the computation 
\begin{equation*}
\int_{0}^{1}h_{\gamma }\left( x\right) dx=\int_{0}^{1}x^{-2\gamma
}dx+\int_{0}^{1}x^{-2\gamma +1}dx=\frac{1}{\left( 1-2\gamma \right) \left(
2-2\gamma \right) }.
\end{equation*}

\newpage

\bibliographystyle{alea3}
\bibliography{example}

\newpage

\section*{Acknowledgements}
The authors thank the anonymous referees for
pointing out the mistakes in the initial version.
The first author gratefully acknowledges the support for this research from
Professor Soledad Torres and the Centro de Investigaci\'{o}n y Modelamiento
de Fen\'{o}menos Aleat\'{o}rios (CIMFAV) of the Universidad de Valpara\'{\i}%
so, Chile, as well as from the Ing\'{e}nieur program at Ecole Polytechnique
in Palaiseau, France. The second author's research was partially supported
by NSF grant DMS 0907321, and by the CIMFAV.

\end{document}